\documentclass[12pt]{amsart}
\usepackage{amssymb}
\usepackage{amsfonts}
\usepackage{amscd}
\usepackage{latexsym}
\usepackage{enumerate}
\usepackage{amsmath}
\usepackage{xy} \xyoption{all}

\numberwithin{equation}{section}

\newcounter{cs}
\stepcounter{cs}
\newcounter{ds}
\stepcounter{ds}

\newcommand{\casos}{\begin{itemize}}
\newcommand{\fcasos}{\end{itemize}\setcounter{cs}{1}}

\newcommand{\Hom}{\mbox{\rm Hom}}

\newcommand{\ol}{\overline}

\newtheorem{lem}{Lemma}[section]
\newtheorem{cor}[lem]{Corollary}
\newtheorem{theor}[lem]{Theorem}
\newtheorem{prop}[lem]{Proposition}

\theoremstyle{definition}
\newtheorem{defi}[lem]{Definition}

\newtheorem{nota}[lem]{Notation}
\newtheorem{remark}[lem]{Remark}

\newtheorem{obs}[lem]{Observation}

\newcommand{\N}{\mathbb{N}}

\newcommand{\Z}{\mathbb{Z}}

     \newcommand{\id}[1][]{\mathbf{1}_{#1}}             
     \newcommand{\mat}[2][]{M_{#1}(#2)}                 
     \newcommand{\gl}[2][]{GL_{#1}(#2)}                 
     \newcommand{\idem}[1]{\mathrm{Idem}(#1)}           
     \newcommand{\expesq}[2]{{\vphantom{#1}}^{#2}\!{#1}}
     \newcommand{\subesq}[2]{{\vphantom{#1}}_{#2}{#1}}  
     \newcommand{\ord}[1]{o(#1)}
     \newcommand{\supp}[1]{\mathrm{supp}(#1)}
     \newcommand{\rang}[2]{\mathrm{rank}_{#1}(#2)}
     \DeclareMathOperator{\diag}{diag}
     \DeclareMathOperator{\Endo}{End}

     \DeclareMathOperator{\coker}{coker}

     \newcommand{\mon}[1]{\mathcal{V}(#1)}              
     \newcommand{\soc}[1]{\mbox{\rm Soc}(#1)}
     \newcommand{\uni}[1]{U(#1)}                        
     \newcommand{\einv}[2]{T(#1\subseteq#2)}            
     \newcommand{\minv}[2]{\Sigma(#1\subseteq#2)}
     \newcommand{\divc}[2]{\mathcal{D}(#1\subseteq#2)}
     \newcommand{\ratc}[2]{\mathcal{R}(#1\subseteq#2)}
     \newcommand{\lpr}[1]{\mathcal{L}_{#1}}             
     \newcommand{\rpr}[1]{\mathcal{R}_{#1}}             
     \newcommand{\cb}[0]{K}                             

     \newcommand{\lanu}[2][]{\ell.\mathrm{ann}_{#1}(#2)}

     \newcommand{\cam}[2][]{P_{#1}(#2)}
     \newcommand{\rat}[1]{P_{\mathrm{\!rat}}(#1)}
     \newcommand{\ser}[1]{P((#1))}
     \newcommand{\aug}[2][]{\varepsilon_{#1}(#2)}
     \newcommand{\tra}[2][e]{\left(#2\right)\delta_{#1}}

\begin{document}
\title[The regular algebra of a quiver]{The regular algebra of a quiver}
\author{Pere Ara and Miquel Brustenga}\address{Departament de Matem\`atiques, Universitat Aut\`onoma de
Barcelona, 08193, Bellaterra (Barcelona),
Spain}\email{para@mat.uab.es, mbrusten@mat.uab.es}

\thanks{Both authors were partially supported by the DGI
and European Regional Development Fund, jointly, through Project
MTM2005-00934, and by the Comissionat per Universitats i Recerca
de la Generalitat de Catalunya.} \subjclass[2000]{Primary 16D70;
Secondary 06A12, 06F05, 46L80} \keywords{von Neumann regular ring,
path algebra, Leavitt path algebra, universal localization}
\date{\today}

\begin{abstract} Let $K$ be a fixed field. We attach to each column-finite
quiver $E$ a von Neumann regular $K$-algebra $Q(E)$ in a
functorial way. The algebra $Q(E)$ is a universal localization of
the usual path algebra $P(E)$ associated with $E$. The monoid of
isomorphism classes of finitely generated projective right
$Q(E)$-modules is explicitly computed.
\end{abstract}
\maketitle

\section*{Introduction}

In a series of recent papers \cite{AA1}, \cite{AA2}, \cite{AMP},
\cite{APS} different aspects of the algebraic structure of the
so-called Leavitt path algebras $L_K(E)$ have been analyzed. These
algebras, defined with coefficients in an arbitrary field $K$, are
the purely algebraic analogues of the important class of
Cuntz-Krieger graph $C^*$-algebras. See the book of Raeburn
\cite{Raeburn} for a recent account about Cuntz-Krieger algebras.

Fix a field $K$. For a column-finite quiver $E$, denote by $P(E)$
the usual path $K$-algebra associated to $E$ and by $L(E)$ the
Leavitt path $K$-algebra of $E$. In this paper, we show that
 the algebra $L(E)$ can be embedded in a von Neumann regular algebra
$Q(E)$ in such a way that the embedding preserves the monoid of
isomorphism classes of finitely generated projective right
modules, that is the inclusion $L(E)\to Q(E)$ induces a monoid
isomorphism  $\mon{L(E)}\cong \mon{Q(E)}$. Moreover the von
Neumann regular algebra $Q(E)$ is obtained from $P(E)$ and also
from $L(E)$ by universal localization, so that it is a
(generalized) ring of fractions of both algebras. Using this we
prove that the construction of $Q(E)$ is functorial with respect
to {\em complete} graph homomorphisms (see Section 4 for the
definition). This enables us to extend many results from the case
of a finite quiver to the case of an arbitrary column-finite
quiver.

Our construction is relevant to the following realization problem
for von Neumann regular rings, which is a variant of a problem
posed by Goodearl in \cite[Fundamental Open Problem]{directsum}.
Recall that an abelian  monoid $M$ is conical in case $x+y=0$
implies $x=y=0$ and that it is a {\em refinement monoid} in case
any equality $x_1+x_2=y_1+y_2$ admits a refinement, that is, there
are $x_{ij}$, $1\le i,j \le 2$ such that $x_i=x_{i1}+x_{i2}$ and
$y_j=x_{1j}+x_{2j}$ for all $i,j$, see e.g. \cite{AMP}.

\medskip

{\bf Realization Problem for von Neumann Regular Rings} Let $M$ be
a countable refinement conical abelian monoid. Is there a von
Neumann regular ring $R$ such that $\mon{R}\cong M$?

\medskip

For every von Neumann regular ring $R$, it is known that $\mon{R}$
is a conical refinement abelian monoid. Indeed this is the case
for the larger class of exchange rings, by \cite[Corollary
1.3]{agop}.

The countability condition is important here. Fred Wehrung
\cite{Wehisrael} proved that there are (even cancellative)
refinement conical abelian monoids of size $\aleph _2$ such that
cannot be realized as $\mon{R}$ for any von Neumann regular ring
$R$. Note that, by the results in \cite{Wehsemforum}, an
affirmative answer to the above realization problem would give a
negative answer to the Fundamental Separativity Problem for von
Neumann regular rings \cite{agop}.

Since the monoids $\mon{L(E)}$ corresponding to column-finite
quivers $E$ were neatly computed in \cite{AMP}, and their
properties are fairly well understood (see \cite[Section 5]{AMP}),
the results in the present paper represent a significant
contribution to the realization problem. The only systematic
approaches to the realization problem known to the authors are the
well-known realization theorem for dimension monoids
(\cite[Theorem 15.24(b)]{vnrr}), and the realization theorem given
in \cite[Theorem 8.4]{AGP}, saying that every countable abelian
group $G$ can be obtained as $K_0(R)$ for some purely infinite
simple regular ring $R$. Since $\mon{R}=\{0\}\sqcup K_0(R)$ for
every purely infinite simple ring $R$, we get that all the monoids
of the form $\{0\}\sqcup G$, for $G$ a countable abelian group,
can be realized as monoids associated to a purely infinite simple
regular ring.

Our results are a generalization of the ones obtained in
\cite{AGP}, where a von Neumann regular envelope of the Leavitt
algebra of type $(1,n)$ was constructed. The Leavitt algebra of
type $(1,n)$ can be seen as the Leavitt path algebra associated
with a quiver with just one vertex and $n+1$ arrows. The
corresponding path algebra is the free associative algebra with
$n+1$ generators, and the construction in \cite{AMP} uses the
algebra of rational formal power series on these generators
(\cite{free},\cite{BR}). A large part of the present paper is
devoted to generalize properties of the (rational) formal power
series algebra to the more general setting of quiver algebras.

We now summarize the content of the paper. In Section 1 we review
some basic concepts and we generalize some results known for the
free algebra to the context of path algebras. In particular we
study the algebra of formal power series on the quiver $E$ and the
algebra $\rat{E}$ of rational series over $E$, which is defined as
the division closure of $\cam{E}$ in $\ser{E}$.

Section 2 contains the basic constructions of our algebras of
Leavitt type, associated with an algebra $R$ which is a subalgebra
of the algebra $\ser{E}$ of formal power series on a finite quiver
$E$ containing the path algebra $\cam{E}$. We show that if the
algebra $R$ is closed under inversion in $\ser{E}$ then the
resulting ring of Leavitt type is von Neumann regular. When we
take $R=P(E)$ we recover the usual Leavitt path algebra $L(E)$
(which is not von Neumann regular in general). When we take
$R=\rat{E}$, the algebra of rational power series on $E$, we get
what we call the {\em regular algebra} $Q(E)$ associated with $E$.

Section 3 contains the computation of the monoid of finitely
generated projective modules over the von Neumann regular algebras
$T$ of Leavitt type constructed in Section 2. In particular we get
that the inclusion $L(E)\to Q(E)$ induces a monoid isomorphism
$\mon{L(E)}\cong \mon{Q(E)}$. Section 4 gives the functoriality of
the construction with respect to complete graph homomorphisms,
which enables us to extend the construction and the computations
from finite to column-finite quivers.

\section{A universal localization of the path algebra of a quiver}

Let $R$ be a ring. We will use the notation $R^n$ (respectively,
$\expesq{R}{n}$) for the left (respectively, right) $R$-module of
$n$-rows (respectively, $n$-columns) with coefficients in $R$. We
will use $M_{m\times n}(R)$ for the space of matrices of size
$m\times n$ over $R$ and $M_n(R)$ for the ring of $n\times n$
matrices over $R$.



In the following, $\cb$ will denote a fixed field and
$E=(E^0,E^1,r,s)$ a finite quiver (oriented graph) with
$E^0=\{1,\dotsc,d\}$. Here $s(e)$ is the {\em source vertex} of
the arrow $e$, and $r(e)$ is the {\em range vertex} of $e$. A {\em
path} in $E$ is either an ordered sequence of arrows
$\alpha=e_1\dotsb e_n$ with $r(e_t)=s(e_{t+1})$ for $1\leqslant
t<n$, or a path of length $0$ corresponding to a vertex $i\in
E^0$, which will be denoted by $p_i$. The paths $p_i$ are called
trivial paths, and we have $r(p_i)=s(p_i)=i$. A non-trivial path
$\alpha=e_1\dotsb e_n$ has length $n$ and we define
$s(\alpha)=s(e_1)$ and $r(\alpha)=r(e_n)$. We will denote the
length of a path $\alpha$ by $|\alpha|$, the set of all paths of
length $n$ by $E^n$, for $n>1$, and the set of all paths by $E^*$.

Let $\cam{E}$ be the $\cb$-vector space with basis $E^*$. It is
easy to see that $\cam{E}$ has a structure of $\cb$-algebra (see
for example \cite[Section III.1]{AuslanderReitenSmalo}), which is
called the {\em path algebra}. Indeed, $\cam{E}$ is the
$\cb$-algebra given by free generators $\{p_i\mid i\in E^0\}\cup
E^1$ and relations:
\begin{enumerate}[(i)]
\item $p_ip_j=\delta_{ij}p_i$ for all $i,\,j\in E^0$. \item
$p_{s(e)}e=ep_{r(e)}=e$ for all $e\in E^1$.
\end{enumerate}




Observe that $A=\oplus_{i\in E^0}\cb p_i\subseteq \cam{E}$ is a
subring isomorphic to $\cb^d$. In general we will identify
$A\subseteq \cam{E}$ with $\cb^d$. An element in $\cam{E}$ can be
written in a unique way as a finite sum $\sum_{\gamma\in
E^*}\lambda_\gamma\gamma$ with $\lambda_\gamma\in\cb$. We will
denote by $\varepsilon$ the augmentation homomorphism, which is
the ring homomorphism $\varepsilon\colon \cam{E}\to\cb^d\subseteq
\cam{E}$ defined by $\aug{\sum_{\gamma\in
E^*}\lambda_\gamma\gamma}=\sum_{\gamma\in
E^0}\lambda_\gamma\gamma$.

\begin{defi} Let $I=\ker (\varepsilon)$ be the augmentation ideal of $\cam{E}$.
Then the {\em $K$-algebra of formal power series} over the quiver
$E$, denoted by $\ser{E}$, is the $I$-adic completion of
$\cam{E}$, that is $\ser{E}\cong \varprojlim \cam{E}/I^n .$
\end{defi}

An element of $\ser{E}$ can be written in a unique way as a
possibly infinite sum $\sum _{\gamma\in E^*}\lambda
_{\gamma}\gamma$ with $\lambda _{\gamma}\in K$. We will also
denote by $\varepsilon $ the augmentation homomorphism on
$\ser{E}$.

Set $R=P(E)$ or $P((E))$.  Observe that the elements in
$\mat[n]{R}$ (or in $R^n$) can also be uniquely written as
possibly infinite sums $\sum_{\gamma\in E^*}\lambda_\gamma\gamma$
with $\lambda_\gamma\in\mat[n]{\cb}$ (respectively with
$\lambda_\gamma\in \cb^n$) and so we can also define over them the
augmentation homomorphisms, which will be denoted also by
$\varepsilon$.

Given an element $r=\sum_{\gamma\in E^*}\lambda_\gamma\gamma$ in
$R$, $\mat[n]{R}$ or $R^n$ we define its support as
$\supp{r}=\{\gamma\in E^*\mid \lambda_\gamma\neq 0\}$ and we
define its order $\ord{r}$ as the minimum length of the paths in
$\supp{r}$.

Define, for $e\in E^1$, the following additive mappings:
$$
\begin{array}{rccc}
\delta_{e}\colon&R&\longrightarrow&R\\
&\displaystyle\sum_{\alpha\in E^*}\lambda_\alpha \alpha&\longmapsto
&\displaystyle\sum_{\substack{\alpha\in E^*\\r(\alpha)=s(e)}}
\lambda_{\alpha e} \alpha
\end{array}
$$
We will write $\delta_{e}$ on the right of its argument. We will
sometimes refer to the maps $\delta_{e}$ as the (left)
transductions. There are corresponding maps on $\mat[n]{R}$ and on
$R^n$, defined componentwise, denoted also by $\delta_{e}$.

The right transductions $\tilde{\delta}_e\colon R\to R$ are
defined similarly by $$\tilde{\delta_e} (\sum_{\alpha\in
E^*}\lambda_\alpha \alpha)=\sum_{\substack{\alpha\in
E^*\\s(\alpha)=r(e)}} \lambda_{e\alpha } \alpha.$$




The following well-known result follows easily from \cite[Theorem
5.3]{Bergman}.
\begin{prop}\label{camEmonoid}
We have an isomorphism $\mon{\cam{E}}\overset{\mon {\varepsilon
}}{\cong}\mon{\cb^d}\cong (\Z ^+)^d$ induced by the augmentation
homomorphism $\varepsilon \colon \cam{E}\longrightarrow \cb^d$.
\end{prop}

For a f.g. projective $\cam{E}$-module $M$, we write
$\rang{\cam{E}}{M}\in (\Z ^+)^d$ for the image of $[M]$ under the
isomorphism $\mon{\cam{E}}\cong (\Z ^+)^d$ of Proposition
\ref{camEmonoid}. Note that $\rang{\cam{E}}{M}=(r_1,\dots ,r_d)$
if and only if $M\cong (p_1\cam{E})^{r_1}\oplus \cdots \oplus
(p_d\cam{E})^{r_d}$. Similarly, we use the notation
$\rang{\cb^d}{M}\in (\Z ^+)^d$ for a f.g. (projective)
$\cb^d$-module $M$.

\begin{defi}
Let $R$ be a subalgebra of $\ser{E}$ closed under the left
transductions $\delta_{e}$ and let $B$ be a right $R$-submodule of
$\expesq{R}{n}$. We say that $B$ is {\em regular} if for every
$b\in B$ with $\ord{b}>0$, we have $\tra{b}\in B$ for all arrows
$e\in E^1$.
\end{defi}
The regularity passes to direct summands:

\begin{lem}\label{lema:submodregular}
Given a regular module $B=B_1\oplus B_2$ we have that $B_1,\,B_2$
are also regular modules.
\end{lem}
\begin{proof}
Given $b\in B_1$ with $\ord{b}>0$, we have $\tra{b}\in B$ for each
arrow $e\in E^1$ by regularity of $B$. So $\tra{b}=b_1^e+b_2^e$
for some unique $b_i^e\in B_i$ with $b^e_ip_{s(e)}=b^e_i$. Thus
$b=\sum_e (b\delta _e )\cdot e=\sum_e b_1^ee+\sum_e b_2^ee$ and we
get that $\sum_e b_2^ee=0$. It follows that $b_2^e=0$ for every
$e\in E^1$ and that $\tra{b}=b_1^e\in B_1$.
\end{proof}

The following result provides a generalization of \cite[Theorem
VII.3.1]{BR}.

\begin{theor}[Characterization of regular modules]
Let $B\subseteq\expesq{\cam{E}}{n}$ be a right $\cam{E}$-module.
Then $B$ is regular if and only if $B=\bigoplus_{i=1}^t B_i$ where
each $B_i$ is a cyclic projective $\cam{E}$-module, $t\le n$, and
$\rang{\cam{E}}{B}=\rang{\cb^d}{\aug{B}}$.
\end{theor}
\begin{proof}
Denote by $\varepsilon_j$ the following composition
$\expesq{\cam{E}}{n}\to\expesq{(\cb^d)}{n}\cong(\expesq{\cb}{n})^d\to\expesq{\cb}{n}$,
where the last homomorphism is the projection onto the $j$-th
component of $(\expesq{\cb}{n})^d$.

We first show that there are $t_1,\dotsc,t_d\in\N$ and
$v_1^{(j)},\dotsc,v_{t_j}^{(j)}\in Bp_j$, where $j=1,2,\dotsc,d$
such that
\begin{enumerate}[(i)]
\item\label{claim1}
$\deg(v_1^{(j)})\leqslant\dotsb\leqslant\deg(v^{(j)}_{t_j}).$
\item\label{claim2} The vectors
$\{\aug[j]{v_i^{(j)}}\}_{i=1,\dotsc,t_j}$ form a $\cb$-basis of
$\aug[j]{Bp_j}$. \item\label{claim3} If $v\in Bp_j$ and
$\deg(v)<\deg(v_i^{(j)})$, then $\aug[j]{v}$ is a linear
combination of
$\aug[j]{v_1^{(j)}},\dotsc,\linebreak\aug[j]{v_{i-1}^{(j)}}$.
\end{enumerate}

For fixed $j$ and $r=-1,0,1,2,\dotsc$ write
$$
F(r)=\{\aug[j]{v}\in\expesq{K}{n}\mid v\in Bp_j,\,\deg(v)\leqslant
r\}.
$$
Observe that $F(r)$  are $\cb$-vector spaces. We have inclusions
$0=F(-1)\subseteq F(0)\subseteq\dotsb\subseteq
F(r)\subseteq\dotsb$. Take integer numbers $0\leqslant
r_1<\dotsb<r_q$ such that $F(r_i-1)\subset F(r_i)$ and that, for
all $r$ with $F(r)\ne 0$, we have $F(r)=F(r_i)$ for some $i$. In
particular, $F(r_q)=\aug[j]{Bp_j}$.

Choose now a $\cb$-basis $\mathcal{B}_1$ of $F(r_1)$, a
$\cb$-basis $\mathcal{B}_2$ of $F(r_2)$ modulo $F(r_1)$ and in
general a $\cb$-basis $\mathcal{B}_i$ of $F(r_i)$ modulo
$F(r_{i-1})$. By definition of the $F$'s, for each $i=1,\dotsc,q$
we can find elements $w_{i,1},\dotsc,w_{i,k_i}\in B$ of degree
less than or equal to $r_i$ such that
$\{\aug[j]{w_{i,1}},\dotsc,\aug[j]{w_{i,k_i}}\}=\mathcal{B}_i$.
Indeed, the degree of each $w_{i,j}$ is exactly $r_i$, otherwise
we would have that $\aug[j]{w_{i,j}}$ would belong to
$F(r_i-1)=F(r_{i-1})$, contradicting that $\mathcal{B}_i$ is a
basis modulo $F(r_{i-1})$. Now put $t_j=k_1+k_2+\dotsb+k_q$ and
define  $v_1^{(j)},\dotsc,v_{t_j}^{(j)}$ by
$$
(v_1^{(j)},\dotsc,v_{t_j}^{(j)})=(w_{1,1},\dotsc,w_{1,k_1},w_{2,1},\dotsc,w_{2,k_2},\dotsc,w_{q,k_q}).
$$
It is clear from this definition that condition \eqref{claim1} is
satisfied. Moreover, since $F(r_q)=\aug[j]{Bp_j}$, we see that
condition \eqref{claim2} is also satisfied. Let $v\in Bp_j$ with
$\deg(v)<\deg(v_i^{(j)})$. Then $v_i^{(j)}=w_{\ell,m}$ for some
$\ell,m$, so that $\deg(v)<r_\ell=\deg(w_{\ell,m})$, and we
conclude that $\aug[j]{v}\in F(r_\ell-1)=F(r_{\ell-1})$ and that
$\aug[j]{v}$ is a $\cb$-linear combination of $\varepsilon
(w_{1,1}),\dotsc,\varepsilon (w_{\ell-1,k_{\ell-1}})$, and so of
$\varepsilon (v_1^{(j)}),\dotsc,\varepsilon (v_{i-1}^{(j)})$. This
shows \eqref{claim3}.

Set $t=\max\{t_1,\dotsc,t_d\}$. Since
$t_j=\dim_K\aug[j]{Bp_j}\leqslant n$ we have that $t\leqslant n$.
For $i=1,\dotsc,t$, take the following submodules of $B$
$$
B_i=v_i^{(1)}\cam{E}+v_i^{(2)}\cam{E}+\dotsb+v_i^{(d)}\cam{E}\subseteq
B
$$
where $v_i^{(j)}=0$ if $i>t_j$. We will now see that $B_i$  are
cyclic projective modules.

For $i=1,\dotsc,t$ and $j=1,\dotsc, d$ set
$$
a_{ij}=\begin{cases}0\text{ if $i>t_j$}\\1\text{ if $i\leqslant
t_j$}\end{cases}
$$
With this notation, for $i=1,\dotsc t$, we have isomorphisms
$B_i\cong a_{i1}(p_1\cam{E})\oplus
a_{i2}(p_2\cam{E})\oplus\dotsb\oplus a_{id}(p_d\cam{E})$. Indeed,
define

\begin{center}
\begin{tabular}{cccc}
$\varphi_i\colon$&$a_{i1}(p_1\cam{E})\oplus\dotsb\oplus a_{id}(p_d\cam{E})$&$\longrightarrow$&$B_i$\\
&$(a_{i1}p_1a_1,\dotsc,a_{id}p_da_d)$&$\longmapsto$&$\displaystyle{\sum_{j=1}^d}v_i^{(j)}p_ja_j$
\end{tabular}
\end{center}
It is clear that the maps $\varphi_i$ are surjective since
$a_{ij}=0$ if and only if $v_i^{(j)}=0$. We claim that
\begin{equation*}\tag{1.1}\label{equ:combi} \sum
_{j=1}^d\sum_{i=1}^{t_j}v_i^{(j)}(p_jb_{ij})=0
\end{equation*}
implies $p_jb_{ij}=0$ for all $i,j$. This will show at once that
each $\varphi _i$ is injective and that the sum $B_1+\cdots +B_t$
is a direct sum.

To prove the claim, we proceed by way of contradiction, so let
\begin{equation*}\tag{1.2}\label{equ:combinalin}
\sum _{j=1}^d\sum_{i=1}^{t_j}v_i^{(j)}(p_jb_{ij})=0\quad\text{with
$0\leqslant\max_{i,j}\{\deg(p_jb_{ij})\}$ smallest possible.}
\end{equation*}

For all $j$ we have
$\sum_{i=1}^{t_j}\aug[j]{v_i^{(j)}}\aug[j]{p_jb_{ij}}=0$ and so we
deduce from (\ref{claim2}) that, for all $i$ and $j$,
$\aug[j]{p_jb_{ij}}=0$ so that $\aug{p_jb_{ij}}=0$, since
$\aug[k]{p_j}=0$ if $k\neq j$. Since these elements have zero
augmentation we can apply the transductions so reducing the
degree:
$$
\sum _{j=1}^d\sum_{i=1}^{t_j}v_i^{(j)}\cdot (\tra{p_jb_{ij}})=0.$$
Observe that $\tra{p_jb_{ij}}=p_j(\tra{p_jb_{ij}})$ for all $e$
and that $\tra{p_jb_{ij}}\ne 0$ for some $i,j$ and some $e\in
E^1$. This leads to a contradiction to the minimality of the
degree in \eqref{equ:combinalin}.

Now set $B'=\bigoplus_{i=1}^tB_i\subseteq B$. We want to see that
$B'=B$. Suppose there is $v\in B\setminus B'$, which we take of
minimal degree. We can write $v=vp_1+\dotsb+vp_d$. Denote
$\deg(v)$ by $d_v$. Let $j$ be an integer such that
$\deg(vp_j)=d_v$ and let $i$ be the smallest integer such that
$\deg(vp_j)<\deg(v_i^{(j)})$ ($i=t_j+1$ if no such integer
exists). By \eqref{claim2} or \eqref{claim3}, depending on  the
case, we have that
$\aug[j]{vp_j}=\sum_{k=1}^{i-1}\lambda_k\aug[j]{v_k^{(j)}}$ and
thus $v'=vp_j-\sum_{k=1}^{i-1}\lambda_kv_k^{(j)}\in B$ satisfies
that $\deg(v')\leqslant\deg(v)$ and $v'\in Bp_j$ with
$\aug[j]{v'}=0$. It follows that $\aug{v'}=0$. Since $B$ is a
regular module we have  $\tra{v'}\in B$ for all $e\in E^1$. By the
minimality of the degree we have that $\tra{v'}\in B'$ and then
$v'=\sum_{e\in E^1}(\tra{v'})e\in B'$. We get from this that
$vp_j\in B'$ for all $j$ such that $\deg(vp_j)=d_v$. Since
$$
v-\sum_{\deg(vp_j)=d_v}vp_j=\sum_{\deg(vp_k)<d_v}vp_k,
$$
we have that
$$
\deg\left(v-\sum_{\deg(vp_j)=d_v}vp_j\right)<d_v
$$
and, again by the minimality of the degree, we get that $v\in B'$.

We now prove the converse. Suppose that $B=\bigoplus_{i=1}^tB_i$
where each $B_i$ is a cyclic projective module, $t\le n$ and
$\rang{\cam{E}}{B}=\rang{\cb^d}{\aug{B}}$. Observe that we have
$$B=\bigoplus _{j=1}^d \bigoplus _{i=1}^{t_j}x_i^{(j)}\cam{E},$$
with $x_i^{(j)}=x_i^{(j)}p_j$ for all $i,j$, where
$\rang{\cam{E}}{B}=(t_1,\dots ,t_d)=\rang{\cb^d}{\varepsilon
(B)}$. It follows from the latter equality that $\{ \varepsilon
_j(x_i^{(j)})\mid i=1,\dots , t_j\}$ is a linearly independent
family of vectors in $\expesq{\cb}{n}$.

Let $v\in B$ such that $\aug{v}=0$. Then for all $j$ we have
$$0=\aug[j]{v}=\sum _{i=1}^{t_j}\aug[j]{x_i^{(j)}}\aug[j]{z_j^i},$$
where $v=\sum _{i,j}x_i^{(j)}z_j^i$, with $z_j^i\in p_j\cam{E}$.
Since $\{\aug[j]{x_i^{(j)}}\mid i=1,\dots ,t_j\}$ are $K$-linearly
independent, we get $\aug[j]{z_j^i}=0$ for all $i,j$, but
$z_j^i=p_jz_j^i$ so $\aug{z_j^i}=0$ and thus $v\delta_e=\sum
_{i,j}x_i^{(j)}\cdot ((z_j^i)\delta _e)\in B$.
\end{proof}

\begin{cor}\label{cor:carac}
If $B\subseteq\expesq{\cam{E}}{n}$ is a regular module, then there
exists $u\in\mat[n]{\cam{E}}$ such that $B=u\expesq{\cam{E}}{n}$.
\end{cor}

The proof of next lemma is standard, see for example \cite[pp.
284--285]{Cohnenc}.

\begin{lem}[Higman's trick]\label{lema:higman}
Given a matrix $M\in\mat[n\times m]{\cam{E}}$, there exist
$\ell\in\N$, $P\in\gl[n+\ell]{\cam{E}}$ and
$Q\in\gl[m+\ell]{\cam{E}}$ such that
$P\left(\begin{smallmatrix}M&0\\0&\id[\ell]\end{smallmatrix}\right)Q$
is a linear matrix.
\end{lem}

\begin{lem}\label{lema:invers}
We have
$\gl[n]{\ser{E}}=\{M\in\mat[n]{\ser{E}}\mid\aug{M}\in\gl[n]{\cb^d}\}$.
\end{lem}
\begin{proof}
If $M\in\gl[n]{\ser{E}}$ then clearly $\aug{M}\in\gl[n]{\cb^d}$.
Let now $M\in\mat[n]{\ser{E}}$ such that
$\aug{M}\in\gl[n]{\cb^d}$. We can write $M=\aug{M}-D$ for some
$D\in\mat[n]{\ser{E}}$ with $\ord{D}>0$. We have that
$$
M^{-1}=\aug{M}^{-1}(\id[n]+D\aug{M}^{-1}+(D\aug{M}^{-1})^2+\dotsb).
$$
\end{proof}

%

\begin{lem}\label{lema:vNreg}
Let $u_0=p-D\in\mat[n]{\cam{E}}$, where
$p\in\idem{\mat[n]{\cb^d}}$ and $D$ is homogeneous of degree $1$.
Suppose further that $B=u_0\expesq{\cam{E}}{n}$ is a regular
$\cam{E}$-module. Then there exist $u\in\mat[n]{\cam{E}}$ and
$v\in\mat[n]{\ser{E}}$ such that $uvu=u$, $vuv=v$ with
$B=u\expesq{\cam{E}}{n}$ and $vu\in\mat[n]{\cam{E}}$.
\end{lem}
\begin{proof}
If the matrix $u_0(1-p)=-D(1-p)$ is nonzero then it is a
homogeneous matrix of degree $1$. The columns of this matrix are
elements of $B$, that is, denoting the elements of the canonical
basis of $\expesq{\cam{E}}{n}$ by $E_i$, we have $u_0(1-p)E_i\in
B$ for all $i=1,\dotsc,n$. These elements are of positive order,
so they decompose as follows:
$$
u_0(1-p)E_i=\sum_{e\in E^1}u_{0e}^ie\quad\text{with unique
$u_{0e}^i\in \expesq{\left(\cam{E}p_{s(e)}\right)}{n}$ ,}
$$
and indeed, being $B$ a regular module, we get $u_{0e}^i\in B$. We
also have that $\deg(u_{0e}^i)<1$, that is, they are elements of
$\aug{B}\cap B$. Since $\aug{B}=p\expesq{(\cb^d)}{n}$ we have
$pu_{0e}^i=u_{0e}^i$ for all $i$ and $e$. In particular,
$(1-p)D(1-p)=0$. Consider the $\cb^d$-submodule  $V_1$ of
$\aug{B}$ generated by $\{u^i_{0e}\mid e\in E^1,\,i=1,\dotsc,n\}$.
Since $\cb^d$ is a semisimple ring, there exists
$q_1\in\idem{\mat[n]{\cb^d}}$, $q_1\leqslant p$, such that
$V_1=q_1\expesq{(\cb^d)}{n}$. Therefore by the above we have that
$q_1\expesq{\cam{E}}{n}\subseteq B$ and, in conclusion we have
seen that
$$
u_0=p-pDp-q_1D(1-p)-(1-p)Dp.
$$

Setting $u_1=(1-q_1)u_0=(p-q_1)-(p-q_1)Dp-(1-p)Dp$, we obtain from
the modular law:
$$B=q_1\expesq{\cam{E}}{n}\oplus u_1\expesq{\cam{E}}{n}.$$

By Lemma \ref{lema:submodregular}, $u_1\expesq{\cam{E}}{n}$ is
again a regular $\cam{E}$-module, so that we can repeat the above
process with  $u_1$. We have that $u_1q_1=-(p-q_1)Dq_1-(1-p)Dq_1$.
As before, for $i=1,\dotsc,n$, we get $u_1q_1E_i=\sum_{e\in
E^1}u^i_{1e}e$ with $u_{1e}^i\in B$ and since we know by degree
considerations that $u_{1e}^i\in\aug{B}$ we obtain that
$(1-p)Dq_1=0$.

Consider the $\cb^d$-submodule $V_2$ of $\aug{B}$ generated by
$\{u^i_{1e}\mid e\in E^1,\,i=1,\dotsc,n\}$. There exists an
idempotent matrix $q_2\in\mat[n]{\cam{E}}$, $q_2\leqslant
(p-q_1)$, such that $V_2=q_2\expesq{(\cb^d)}{n}\leqslant\aug{B}$.
Similarly to the above argument we get
$q_2\expesq{\cam{E}}{n}\subseteq B$ since $u_{1e}^i\in B$. Putting
\begin{align*}
u_2&=(1-(q_1+q_2))u_1=(1-(q_1+q_2))u_0\\
&=(p-(q_1+q_2))-(p-(q_1+q_2))D(p-q_1)-(1-p)D(p-q_1)
\end{align*}
we have that $B=(q_1+q_2)\expesq{\cam{E}}{n}\oplus
u_2\expesq{\cam{E}}{n}$.

Iterating this process we get a sequence of idempotent, pairwise
orthogonal matrices $q_1,\dotsc,q_\ell\in \idem{\mat[n]{K^d}}$,
with $q_i\leqslant p$ for all $i$ in such a way that
$q_i\expesq{\cam{E}}{n}\subseteq B$; we also have independent
$\cb^d$-modules $V_i=q_i\expesq{(\cb^d)}{n}$ and
$u_1,\dotsc,u_\ell\in\mat[n]{\cam{E}}$,
\begin{multline*}
u_i=(p-(q_1+\dotsb+q_i))-(p-(q_1+\dotsb+q_i))D(p-(q_1+\dotsb+q_{i-1}))\\
-(1-p)D(p-(q_1+\dotsb+q_{i-1}))
\end{multline*}
such that $B=(q_1+\cdots +q_i)\expesq{\cam{E}}{n}\oplus
u_i\expesq{\cam{E}}{n}$ for all $i$.

Since we have an ascending chain of submodules of a Noetherian
module:
$$
V_1\subseteq V_1\oplus V_2\subseteq\dotsb\subseteq V_1\oplus\dotsb\oplus V_i\subseteq\dotsb\subseteq\aug{B}
$$
we see that the above process will stop in a finite number of
steps, say after $\ell$ steps. Write $q=q_1+\dotsb+q_\ell$. In
principle we have
$$
u_\ell=(p-q)-(p-q)D(p-(q_1+\dotsb+q_{\ell-1}))-(1-p)D(p-(q_1+\dotsb+q_{\ell-1})),
$$
but since this process stops in the step $\ell$ necessarily we
have $u_\ell q_\ell=0$, that is,
$$
u_\ell=(p-q)-(p-q)D(p-q)-(1-p)D(p-q)
$$
with $B=q\expesq{\cam{E}}{n}\oplus u_\ell\expesq{\cam{E}}{n}$. Now
set
$$
u=q+u_\ell=p-(p-q)D(p-q)-(1-p)D(p-q),
$$
and note that $B=u\expesq{\cam{E}}{n}$. If we set
$$
v=p+(p-q)D(p-q)+\left((p-q)D(p-q)\right)^2+\dotsb\in\mat[n]{\ser{E}}
$$
then $vu=p\in\mat[n]{\cam{E}}$, $uvu=u$ and $vuv=v$.
\end{proof}

\begin{theor}[Stable inertia]\label{tma:inercia}
Let $\subesq{A}{\cam{E}}\subseteq\ser{E}^n$ and
$B_{\cam{E}}\subseteq\expesq{\ser{E}}{n}$ be left and right
$\cam{E}$-submodules, respectively, such that for every $a\in A$
and for every $b\in B$ we have $ab\in\cam{E}$. Then there exist
$m\in\N$ and $u,\,v\in\mat[n+m]{\ser{E}}$ such that for all $a\in
A\oplus\cam{E}^m$ and all $b\in B\oplus\expesq{\cam{E}}{m}$ we
have
$$
au\in\cam{E}^{n+m},\quad vb\in\expesq{\cam{E}}{n+m}\text{ and
}ab=(au)(vb).
$$
\end{theor}
\begin{proof}
We will show first the case where $B\subseteq\expesq{\cam{E}}{n}$.
We can assume that
$$
B=\{b\in\expesq{\cam{E}}{n}\mid\forall a\in A,\,ab\in\cam{E}\}.
$$
Under this assumption we have that $B$ is a regular right
$\cam{E}$-module. Indeed, if we take $b\in B$ with $\ord{b}>0$ we
know that $b=\sum_{e\in E^1}b_ee$ for some unique elements
$b_e\in\expesq{(\cam{E}p_{s(e)})}{n}$. For every $a\in A$ we have
that
$$
ab=a\left(\sum_{e\in E^1} b_ee\right)=\sum_{e\in E^1} (ab_e)e\in\cam{E},
$$
so that $ab_e\in\cam{E}$ for all $e\in E^1$ and for all $a\in A$,
which tells us that $b_e\in B$ for all $e\in E^1$.

By Corollary \ref{cor:carac}, there exists
$u_0\in\mat[n]{\cam{E}}$ such that $B=u_0\expesq{\cam{E}}{n}$. We
get $\aug{B}=\aug{u_0}\expesq{(\cb^d)}{n}$.

By Higman's trick (Lemma \ref{lema:higman}) there exist $m\in\N$
and $P,\,Q\in\gl[n+m]{\cam{E}}$ such that the matrix
$$
u_1=P\begin{pmatrix}u_0&0\\0&\id[m]\end{pmatrix}Q,
$$
is a linear matrix.
Now we consider the left $P(E)$-submodule
$A'=(A\oplus\cam{E}^m)P^{-1}$ of $P((E))^{n+m}$ and the right
$P(E)$-submodule
$B'=u_1\expesq{\cam{E}}{n+m}=P(B\oplus\expesq{\cam{E}}{m})$ of
$\expesq{\cam{E}}{n+m}$, and observe that $B'$ is regular and that
the generating matrix $u_1$ of $B'$ is linear.

Since $\mat[n+m]{\cb^d}$ is unit-regular there exists
$x\in\gl[n+m]{\cb^d}$ such that $\aug{u_1}x\aug{u_1}=\aug{u_1}$.
Thus the idempotent matrix
$p:=\aug{u_1}x\in\idem{\mat[n+m]{\cb^d}}$ satisfies
$\aug{B}=p\expesq{(\cb^d)}{n+m}$. Using $u_1x$ instead of $u_1$ as
a generator of $B'$, we can assume that $u_1=p-D$ with
$p\in\idem{\mat[n+m]{\cb^d}}$ and $D\in\mat[n+m]{\cam{E}}$
homogeneous of degree $1$.

We are now in the hypothesis of Lemma \ref{lema:vNreg}, so that
there exist $u_2\in\mat[n+m]{\cam{E}}$ and
$v_0\in\mat[n+m]{\ser{E}}$ such that $u_2v_0u_2=u_2$,
$v_0u_2v_0=v_0$ with $B'=u_2\expesq{\cam{E}}{n+m}$ and
$v_0u_2\in\mat[n+m]{\cam{E}}$.

Set $u=P^{-1}u_2$ and $v=v_0P$. Given $a\in A\oplus\cam{E}^m$ and
$b\in B\oplus\expesq{\cam{E}}{m}$ we have that $aP^{-1}\in A'$ and
$Pb\in B'$. Thus $Pb=u_2b'$ for some $b'\in\expesq{\cam{E}}{n+m}$.
The following identities hold:
\begin{align*}
(au)(vb)&=((aP^{-1})u_2)(v_0(Pb))=(aP^{-1})((u_2v_0u_2)b')\\
        &=(aP^{-1})(u_2b')=(aP^{-1})(Pb)=ab,\\
vb&=v_0(Pb)=(v_0u_2)b'\in\expesq{\cam{E}}{n+m}.
\end{align*}
Moreover, using that $u_2=u_1y$ for some $y\in\mat[n+m]{\cam{E}}$
we have that
$$
au=(aP^{-1})u_2=(aP^{-1})u_1y=a\begin{pmatrix}u_0&0\\0&\id[m]\end{pmatrix}Qy\in\cam{E}^{n+m}.
$$
We have shown the result in the case where
$B\subseteq\expesq{\cam{E}}{n}$. \vspace{14pt}

Now we shall see that the general case can be reduced to the case
considered before. Taking, if needed, a bigger subset $A$, we can
assume that
$$
A=\{a\in\ser{E}^n\mid\forall b\in B,\,ab\in\cam{E}\}
$$
so that $A$ is a regular left $\cam{E}$-module.

For $i\in\{1,\dotsc,d\}$, put
$$
J_i=\{j\in\{1,\dotsc,n\}\mid\forall b=(b_1,\dotsc,b_n)^t\in B,\,
p_ib_j\in\cam{E}\}.
$$
These sets give us a measure of how far we are from the preceding
situation. Assume that for all $i$ with $1\leqslant i\leqslant d$,
for all $a=(a_1,\dotsc,a_n)\in A$ and for all $j\notin J_i$ we
have $a_jp_i=0$. In this case, we consider the following diagonal
matrix:
$$
q=\diag(d_1,\dotsc,d_n) \text{ where } d_j=\sum_{i\in\{i\mid j\in
J_i\}}p_i.
$$
For all $a=(a_1,\dotsc,a_n)\in A$ and all $b=(b_1,\dotsc,b_n)^t\in
B$  we have
$$
ab=\sum_ia_ib_i=\sum_ia_id_ib_i+\sum_ia_i(1-d_i)b_i=\sum_ia_id_ib_i=a(qb)
$$
with $qb\in\expesq{\cam{E}}{n}$ so that considering $qB$ instead
of $B$ we can reduce ourselves to the above case.

Otherwise, suppose that there exist $i_0\in\{1,\dotsc,d\}$,
$a=(a_1,\dotsc,a_n)\in A$ and $j_0\notin J_{i_0}$ such that
$a_{j_0}p_{i_0}\neq 0$. In particular, we can take
$w\in\supp{a_{j_0}}$ such that $wp_{i_0}\neq 0$. Write
$w={e_1}\dotsm {e_m}$ for some arrows $e_1,\dotsc,e_m\in E^1$ with
$r(e_m)=i_0$. Consider the following set:
$$
S_a:=\bigcup_{i=1}^d\left(\left(\bigcup_{j\notin
J_i}\supp{\aug{a_j}}\right)\cap\{p_i\}\right)
$$

Observe that, if we denote the elements of the canonical basis of
$\cam{E}^n$ by $E_j$, by the definition of $J_i$ we have that
$p_iE_j\in A$ if and only if $j\in J_i$. Thus in case that
$S_a=\varnothing$ we get that $\aug{a}\in A$ and we can assume
that $\ord{a}>0$. Consider the subword $w_1={e_1}\dotsm {e_{m_1}}$
of $w$, where $m_1=\ord{a}$. By the regularity of $A$ we have that
$a'=\tilde{\delta}_{e_{m_1}}\cdots
\tilde{\delta}_{e_2}\tilde{\delta}_{e_1}(a)$ is a nonzero element
in $A$ with $w'={e_{m_1+1}}\dotsm {e_m}\in\supp{a'_{j_0}}$.

Now we can repeat the same argument with $a'$ and $w'$. Since
$wp_{i_0}\neq 0$ and the above process decrease the length of the
word, we will arrive at some $a^{(k)}\in A$ such that
$S_{a^{(k)}}\neq\varnothing$ and therefore we will always be able
to reduce to this case.

If $S_a\neq\varnothing$ we can take $p_{\ell}\in S_a$ in such a
way that for some $j_0\notin J_\ell$ we have
$p_\ell\in\supp{\aug{a_{j_0}}}$. Set
$$
a''
=(p_\ell a_1,\dotsc,p_\ell a_{j_0-1},p_\ell a_{j_0}+p_1+\dotsb+p_{\ell-1}+p_{\ell+1}
+\dotsb+p_r,p_\ell a_{j_0+1},\dotsc,p_\ell a_n).
$$
By Lemma \ref{lema:invers} we have that the $j_0$-th component of
$a''$ is invertible in $\ser{E}$, so that we can consider the
following invertible matrix:
$$
M=
\begin{pmatrix}
1&0&\hdotsfor{2}&0\\
&\ddots&&&\\
a''_1&\hdotsfor{1}&a''_{j_0}&\hdotsfor{1}&a''_n\\
&&&\ddots&\\
0&\hdotsfor{2}&0&1
\end{pmatrix}\in\gl[n]{\ser{E}}
$$
obtained by substituting the $j_0$-th row of the identity matrix
by $a''$.

For all $b=(b_1,\dotsc,b_n)^t\in B$ we have that
$$
Mb=(b_1,\dotsc,b_{j_0-1},b_{j_0}',b_{j_0+1},\dotsc,b_n)^t,$$ where
$b_{j_0}'=p_\ell
ab+(p_1+\dotsb+p_{j_0-1}+p_{j_0+1}+\dotsb+p_n)b_{j_0}$. Hence, if
we substitute $A$ by $AM^{-1}$ and $B$ by $MB$, this will not
change the sets $J_i$ for $i\neq\ell$ but will increase $|J_\ell|$
by one (we now have $j_0\in J_\ell$). Repeating the above process
a finite number of times, we will arrive at the case where
$B\subseteq\expesq{\cam{E}}{n}$.
\end{proof}

Let $R$ be a ring. Recall that an $R$-module $P$ is {\em stably
free} if $P\oplus R^m\cong R^n$ for some  $m,n\in\N$.
\begin{defi} (\cite[page 15]{free})
A ring $R$ is called a {\em Hermite ring} if it has IBN (invariant
basis number) and the stably free modules are free.
\end{defi}

Recall that a {\em full matrix} over a ring $R$ is a square matrix
$A$ over $R$, of size $n\times n$ say, such that $A$ cannot be
written as a product $A=BC$, where $B\in M_{n\times (n-1)}(R)$ and
$C\in M_{(n-1)\times n}(R)$, see \cite[page 159]{free}. We need
the following lemma:

\begin{lem}\label{lema:estableplena}
If $A$ is a full matrix over $P(E)$ then $$A\oplus \id[m] :=
\begin{pmatrix}A&0\\0&\id[m]\end{pmatrix}$$
is also full, for every $m\ge 0$.
\end{lem}

\begin{proof}
It follows from Proposition \ref{camEmonoid} that $P(E)$ is a
Hermite ring. Now the result is a consequence of \cite[Proposition
5.6.2]{free}.
\end{proof}

\begin{defi}(\cite[page 250]{free})
\label{defi:honest} A ring homomorphism $f\colon R\to S$ is {\em
honest} if it sends full matrices over $R$ to full matrices over
$S$.
\end{defi}

\begin{cor}\label{cor:honest}
The inclusion $\cam{E}\hookrightarrow\ser{E}$ is a honest
inclusion.
\end{cor}
\begin{proof}
Let $C\in\mat[n]{\cam{E}}$ be a matrix such that $C$ is not a full
matrix over $\ser{E}$. Then we can write $C=AB$ where
$A\in\mat[n\times\ell]{\ser{E}}$ and $B\in\mat[\ell\times
n]{\ser{E}}$ with $\ell<n$. By Theorem \ref{tma:inercia} there
exist $m\in\N$ and $u,\,v\in\mat[\ell+m]{\ser{E}}$ such that
$$
\begin{pmatrix}C&0\\0&\id[m]\end{pmatrix}=\begin{pmatrix}A&0\\0&\id[m]\end{pmatrix}
\begin{pmatrix}B&0\\0&\id[m]\end{pmatrix}=\left(\begin{pmatrix}A&0\\0&\id[m]
\end{pmatrix}u\right)\left(v\begin{pmatrix}B&0\\0&\id[m]\end{pmatrix}\right),
$$
so that we obtain a decomposition of the matrix $C\oplus\id[m]$ in
$\cam{E}$, showing that it is not full over $\cam{E}$. By Lemma
\ref{lema:estableplena} the matrix $C$ is not full over $\cam{E}$,
as required.
\end{proof}

We recall the following notation and definitions; see for example
\cite[10.2.2]{Luck}.

\begin{nota}
Given a ring $S$ and a subring $R\subseteq S$, we denote by
$\einv{R}{S}$ the set of all the elements of $R$ which are
invertible in $S$, and we denote by $\minv{R}{S}$ the set of all
square matrices over $R$ which are invertible over $S$.
\end{nota}

\begin{defi}
Let $S$ be a ring.
\begin{enumerate}[(i)]
\item A subring $R\subseteq S$ is {\em closed under inversion} in $S$ if
$\einv{R}{S}=\uni{R}$, that is, if $\uni{R}=R\cap\uni{S}$.
\item A subring $R\subseteq S$ is {\em rationally closed} in $S$ if
$\minv{R}{S}=\gl{R}$, that is, if $\gl{R}=\mat{R}\cap\gl{S}$.
\item Given a subring $R\subseteq S$ the {\em division closure} of $R$ in $S$,
denoted by $\divc{R}{S}$, is the smallest subring of $S$ which is
closed under inversion and contains $R$.
\item Given a subring $R\subseteq S$ the {\em rational closure} of $R$ in $S$,
denoted by $\ratc{R}{S}$, is the smallest subring of $S$ which is
rationally closed and contains $R$.
\end{enumerate}
\end{defi}

\begin{defi}
The {\em $\cb$-algebra of rational series over the quiver} $E$,
denoted by $\rat{E}$, is the division closure of $\cam{E}$ in
$\ser{E}$.
\end{defi}

\begin{obs}
\label{rateqdiv} We note that $\rat{E}$ is also the rational
closure of $\cam{E}$. Indeed if $M$ is a matrix over $\rat{E}$
that becomes invertible in $\ser{E}$, then $\varepsilon (M)$ is
invertible over $K^d$, and replacing $M$ with $\varepsilon
(M)^{-1}M$, we may assume that $\varepsilon (M)$ is an identity
matrix. Hence the diagonal entries of $M$ are invertible in
$\ser{E}$ by Lemma \ref{lema:invers} and so are invertible in
$\rat{E}$. By applying to $M$ a suitable sequence of elementary
row transformations, we may further assume that $M$ is diagonal.
It follows that $M$ is invertible over $\rat{E}$, as claimed.
\end{obs}

\begin{obs}
\label{rateqdiv2} The same argument as above shows that if $R$ is
a subalgebra of $\ser{E}$ closed under inversion and containing
$K^d$ then a matrix $M$ over $R$ is invertible over $R$ if and
only if $\varepsilon (M)$ is invertible over $K^d$.
\end{obs}

Set $\Sigma:= \minv{\cam{E}}{\ser{E}} $, and let $\iota \colon
\cam{E}\to \Sigma^{-1}\cam{E}$ be the {\it universal localization}
of $\cam{E}$ with respect to $\Sigma$, cf. \cite{free},
\cite{scho}. By the universal property, we get a unique
$K$-algebra homomorphism $f\colon \Sigma^{-1}\cam{E}\to \rat{E}$
such that $\phi =f\circ \iota$, where we denote by $\phi \colon
\cam{E}\to \rat{E}$ the natural inclusion. It follows from a
well-known general fact (see for instance \cite[Lemma
10.35(3)]{Luck}) that the map $f$ is surjective.

\begin{theor}
\label{ratuniversal} Let $\Sigma=\minv{\cam{E}}{\ser{E}}$. Then
$\rat{E}$ coincides with the universal localization of $\cam{E}$
with respect to $\Sigma$.
\end{theor}
\begin{proof}
Consider $\iota\colon\cam{E}\to\Sigma^{-1}\cam{E}$, the universal
localization of the path algebra with respect to the set $\Sigma$.
As observed above, we have a surjective $K$-algebra homomorphism
$f\colon \Sigma^{-1}\cam{E}\to \rat{E}$, and we want to see that
it is injective.

By Corollary \ref{cor:honest} we have that the inclusion
$\cam{E}\hookrightarrow\ser{E}$ is honest. Moreover, it is
$\Sigma$-inverting  and it is easily seen that $\Sigma $ is
multiplicative and factor-closed (see \cite[Chapter 7]{free} for
the definitions of these concepts). It follows from
\cite[Proposition 7.5.7(ii)]{free} that $f$ is injective.
\end{proof}

\section{Construction of the algebras}

In this section, we will give the basic construction of the
algebras associated with a finite quiver.

\begin{defi}
Given a quiver $E=(E^0,E^1,r,s)$, consider the sets
$\overline{E}^0=E^0$, $\overline{E}^1=\{\ol{e}\mid e\in E^1\}$ and
the maps $\overline{r},\overline{s}\colon \overline{E}^1\to
\overline{E}^0$ defined via $\overline{r}(\ol{e})=s(e)$ and
$\overline{s}(\ol{e})=r(e)$. Define the {\em inverse quiver} of
$E$ as the quiver
$\overline{E}=(\overline{E}^0,\overline{E}^1,\overline{r},\overline{s})$.
\end{defi}

\begin{nota}
Given a path $\alpha=e_1\dotsm e_n\in E^*$ denote by
$\ol{\alpha}=\ol{e}_n\dotsm \ol{e}_1$ the corresponding path in
the inverse quiver. Of course, if $i\in E^0$, then $\ol{p}_i=p_i$.
\end{nota}

Set $R=\cam{E}$ or $\ser{E}$. For $e\in{E}^1$ we define the
following  $\cb$-algebra endomorphism,
$$
\begin{array}{rcccl}
\tau_{e}\colon&R&\longrightarrow&R\\
&p_{s(e)}&\longmapsto &p_{r(e)}\\
&p_{r(e)}&\longmapsto &p_{s(e)}\\
&p_i&\longmapsto &p_i&\quad i\neq s(e),r(e)\\
&{f}&\longmapsto &0 &\quad\forall f\in E^1 \, .
\end{array}
$$
It is clear that they are $\cb$-algebra endomorphisms, since they
are defined by the composition of the augmentation with an
automorphism of $\aug{R}$ and the inclusion of $\aug{R}$ in $R$.
We will write $\tau_{e}$ on the right of its argument (and
compositions will act accordingly).
\begin{defi}
Let $R$ be a ring and $\tau\colon R\to R$ a ring endomorphism. A
left {\em $\tau$-derivation} is an additive mapping $\delta\colon
R\to R$ satisfying  $(rs)\delta=(r\delta)\cdot (s\tau)+r\cdot
(s\delta)$ for all $r,s\in R$.
\end{defi}
\begin{lem}
For every $e\in {E}^1$, $\delta_{e}$ is a left
$\tau_{e}$-derivation.
\end{lem}
\begin{proof}
Set $r=\sum_{\alpha\in E^*}\lambda_\alpha \alpha$ and
$s=\sum_{\beta\in E^*}\mu_\beta \beta$. Its product is
$rs=\sum_{\gamma\in E^*}\nu_\gamma \gamma$ where
$\nu_\gamma=\sum_{\gamma=\alpha\beta}\lambda_\alpha \mu_\beta$. On
one hand we have that, if $s(e)\neq r(e)$,
\begin{align*}
(r\delta _e)\cdot (s\tau _e )&=\left(\sum_{\substack{\alpha\in
E^*\\r(\alpha)=s(e)}}\lambda_{\alpha e}\alpha\right)
\left(\mu_{r(e)}p_{s(e)}+\mu_{s(e)}p_{r(e)}+\sum_{\substack{i\in E^0\\i\neq r(e),s(e)}}\mu_ip_i\right)\\
&=\sum_{\substack{\alpha \in E^*\\r(\alpha)=s(e)}}\left(\lambda_{\alpha e}\mu_{r(e)}\right)\alpha
\end{align*}
and note that, in case  $s(e)=r(e)$, we get indeed the same
expression. Also,
$$
r\cdot (s\delta _e )=\left(\sum_{\alpha\in E^*}\lambda_\alpha
\alpha\right) \left(\sum_{\substack{\beta\in
E^*\\r(\beta)=s(e)}}\mu_{\beta
e}{\beta}\right) =\sum_{\substack{\gamma\in E^*\\
r(\gamma)=s(e)}}\left(\sum_{\gamma=\alpha\beta}
\lambda_\alpha\mu_{\beta e}\right){\gamma}.
$$
On the other hand, we see that
\begin{align*}
\tra{rs}&=\tra{\sum_{\gamma\in E^*}\nu_\gamma \gamma}
=\sum_{\substack{\gamma\in E^*\\r(\gamma)=s(e)}}\nu_{\gamma e} \gamma
=\sum_{\substack{\gamma\in E^*\\r(\gamma)=s(e)}}\left(\sum_{\gamma e
=\alpha\beta}\lambda_\alpha\mu_{\beta}\right)\gamma\\
&=\sum_{\substack{\gamma\in E^*\\r(\gamma)=s(e)}}
\left(\sum_{\gamma=\alpha\beta}\lambda_\alpha\mu_{\beta e}\right){\gamma}
+\sum_{\substack{\gamma\in E^*\\r(\gamma)=s(e)}}(\lambda_{\gamma e}\mu_{r(e)})\gamma.
\end{align*}
Therefore, $\tra{rs}=(r\delta _e)\cdot (s\tau _e )+r\cdot (s\delta
_e)$.
\end{proof}

In the rest of this section, $E$ will denote a finite quiver with
$E^0=\{1,\dots ,d \}$.

\begin{prop}\label{prop:anellE}
Given a quiver $E$ and a $\cb$-subalgebra $R$ of $\ser{E}$,
containing $\cam{E}$ and closed under all the left transductions
$\delta_{e}$, there exists a ring $S$ such that:
\begin{enumerate}[(i)]
\item There are embeddings
\begin{align*}
\begin{aligned}
\mathcal{R}\colon &  R\to S \qquad \text{\rm and}\qquad  z\colon
& P(\overline{E})\to S \\
& r\mapsto \rpr{r}  &\ol{\alpha}\mapsto z_{\ol{\alpha}}
\end{aligned}
\end{align*}
such that $z_{p_i}=\rpr{p_i}$ for all $i$ and
\begin{equation*}\label{eq:skewprod}\tag{2.1}
\rpr{r}\cdot z_{\ol{e}}=z_{\ol{e}}\cdot \rpr{(r\tau
_e)}+\rpr{(r\delta _e)}
\end{equation*}
for all $e\in E^1$ and all $r\in R$. \label{item:prod2}
\item $S$ is projective as a right $R$-module. Indeed,
$S=\oplus_{\gamma\in E^*}S_\gamma$ with $S_\gamma\cong
p_{s(\gamma)}R$ as $R$-modules.\label{item:proj} Moreover, every
element of $S$ can be uniquely written as a finite sum $\sum
_{\gamma\in E^*}z_{\ol{\gamma}}\rpr{a_{\gamma}}$, where
$a_{\gamma}\in p_{s(\gamma)}R$ for all $\gamma \in E^*$.
\end{enumerate}
\end{prop}
\begin{proof}
Set $T=\Endo_{\cb}(R)$. The elements of $T$ will act on the right
of their arguments. For $r\in R$ denote by $\rpr{r}$ the operator
in $T$ given by right multiplication by $r$. The map
$\mathcal{R}\colon R\to T$ is clearly an injective $K$-algebra
morphism.

For each $\ol{e}\in\overline{E}^1$ consider the elements
$z_{\ol{e}}\in T$ defined by $$ (r)z_{\ol{e}}=(r)\delta _e.
$$ Let $S$ be the subring of $T$ generated by $R$ and by all the
elements $z_{\ol{e}}$ defined above. For $\ol{e}\in
\overline{E}^1$, we have
$$z_{\ol{e}}=z_{\ol{e}}\rpr{p_{s(e)}}=\rpr{p_{r(e)}}z_{\ol{e}},$$ so that
there exists a unique $K$-algebra morphism $z\colon
\cam{\ol{E}}\to S$ such that $z(\ol{e})=z_{\ol{e}}$ for all $e\in
E^1$ and $z(p_i)=\rpr{p_i}$ for all $i\in E^0$.

For $r,s\in R$ and $e\in E^1$, we have
\begin{equation*}
\begin{split}
(s)(\rpr{r}z_{\ol{e}}) & =(sr)z_{\ol{e}}=(sr)\delta _e= (s\delta
_e)\cdot
(r\tau _e)+s\cdot (r\delta _e) \\
& = (s)[z_{\ol{e}}\rpr{(r\tau_e)}+\rpr{(r\delta _e)}].
\end{split}
\end{equation*}
Hence,
\begin{equation*}\notag
\rpr{r}\cdot z_{\ol{e}}=z_{\ol{e}}\cdot \rpr{(r\tau
_e)}+\rpr{(r\delta _e)}
\end{equation*}
for all $\ol{e}\in\overline{E}^1$ and all $r\in R$, which shows
the formula \eqref{eq:skewprod}. From this we conclude that $S$ is
generated as a right $R$-module by monomials $z_{\ol{\gamma}}$
where $\gamma\in E^*$, and thus every element in $S$ can be
written as a finite sum $\sum _{\gamma\in E^*}z_{\ol{\gamma}
}\rpr{a_{\gamma}}$, where $a_{\gamma}\in p_{s(\gamma)}R$ for all
$\gamma \in E^*$. It remains to check uniqueness of the
expression. For, assume that we have $\sum _{\gamma\in
E^*}z_{\ol{\gamma}}\rpr{a_{\gamma}}=0$, where $a_{\gamma}\in
p_{s(\gamma)}R$ are not all $0$. Let $\gamma _0\in E^*$ be a path
of minimal length in the support of this expression, so that
$a_{\gamma_0}\ne 0$. Observe that
$$0=(\gamma _0)(\sum _{\gamma\in
E^*}z_{\ol{\gamma}}\rpr{a_{\gamma}})=p_{s(\gamma_0)}a_{\gamma
_0}=a_{\gamma_0},$$ which gives a contradiction. It follows that
every element in $S$ can be uniquely written as a finite sum $\sum
_{\gamma\in E^*}z_{\ol{\gamma}}\rpr{a_{\gamma}}$, where
$a_{\gamma}\in p_{s(\gamma)}R$ for all $\gamma \in E^*$, which
gives (ii) and also gives the injectivity of the map $z\colon
\cam{\ol{E}}\to S$. This completes the proof.
\end{proof}

\begin{nota}
We will denote the ring $S$ of Proposition \ref{prop:anellE} by
$R\left<\ol{E};\tau,\delta\right>$ where $\tau$ and $\delta$ stand
for $(\tau_{e})_{e\in E^1}$ and $(\delta_{e})_{e\in E^1}$,
respectively. Moreover, since the maps $\mathcal{R}\colon R\to S$
and $z\colon P(\ol{E})\to S $ are injective, we will identify the
elements of $R$ and of $P(\ol{E})$ with their images in $S$ under
these maps. Note that the fundamental relation \eqref{eq:skewprod}
in Proposition \ref{prop:anellE} becomes $r\ol{e}=\ol{e}(r\tau
_e)+(r\delta _e)$ for all $e\in E^1$. In particular
$f\ol{e}=\delta_{e,f}p_{s(e)}$ for $e,f\in E^1$. With this
notation, Proposition \ref{prop:anellE}(ii) says that each element
$s\in R\left<\ol{E};\tau,\delta\right>$ can be uniquely written as
a finite sum
\begin{equation*}\label{eq:uniqueform}\tag{2.2}
s=\sum _{\gamma\in E^*}\ol{\gamma}a_{\gamma}\, ,
\end{equation*}
where $a_{\gamma}\in p_{s(\gamma)}R$ for all $\gamma \in E^*$.
\end{nota}

The algebra $R\left<\ol{E};\tau,\delta\right>$ is characterized by
the following universal property:

\begin{prop}\label{prop:panellE}
Let $\phi\colon R\to B$ be a homomorphism of $\cb$-algebras and
assume that for each $e\in E^1$ there exists $t_{\ol{e}}\in
(p_{r(e)}\phi )B (p_{s(e)}\phi )$ such that
$(r\phi)t_{\ol{e}}=t_{\ol{e}}(r\tau_e\phi)+(r\delta _e\phi)$ for
all $r\in R$ and all $e\in E^1$. Then $\phi$ can be uniquely
extended to a $\cb$-algebra homomorphism $\overline{\phi}\colon
R\left<\ol{E};\tau,\delta\right>\to B$ such that
$\ol{e}\overline{\phi}=t_{\ol{e}}$ for all $e\in E^1$.
\end{prop}
\begin{proof}
The proof is similar to the one of \cite[Proposition 3.3]{AGP}.
Set $S=R\left<\ol{E};\tau,\delta\right>$. It is enough to build a
$K$-algebra with the universal property and to show that it is
isomorphic with $S$. Let $F:=R*_{K^d}P(\ol{E})$ be the coproduct
of $R$ and $P(\ol{E})$ over $K^d$, with canonical maps
$\psi_1\colon R\to F$ and $\psi _2\colon P(\ol{E})\to F$.  Let
$S_1$ be the $K$-algebra obtained by imposing the relations
$(r\psi _1)(\ol{e}\psi _2)=(\ol{e}\psi _2)(r\tau _e\psi
_1)+(r\delta _e\psi_1)$ for all $e\in E^1$ and all $r\in R$. Then
clearly $S_1$ satisfies the required universal property. In
particular we get a $K$-algebra homomorphism $S_1\to S$ extending
the canonical maps $R\to S$ and $P(\ol{E})\to S$. Using the
defining relations, we see that every element in $S_1$ can be
written in the form $\sum _{\gamma\in E^*}(\ol{\gamma}\psi _2
)((p_{s(\gamma)}a_{\gamma})\psi _1)$, where $a_{\gamma}\in R$. Now
it follows from Proposition \ref{prop:anellE} that the map $S_1\to
S$ is an isomorphism.\end{proof}

In the following $R$ will denote a $\cb$-subalgebra of $\ser{E}$
containing $\cam{E}$, closed under inversion (in $\ser{E}$) and
closed under all the left transductions $\delta_{e}$. Examples
include the power series algebra $\ser{E}$ and the algebra
$\rat{E}$ of rational series. Indeed this follows from the fact
that elements $a\in \Sigma^{-1}\cam{E}$ can be written as
$$a=bA^{-1}c,$$
where $A\in \Sigma$ and $b$ and $c$ are a row vector and a column
vector of suitable size. This implies that elements in $\rat{E}$
also have this expression (see Theorem \ref{ratuniversal}).
Observe that, for $A\in \Sigma$ we have $(A^{-1})\delta _e=
-A^{-1}\cdot (A\delta _e)\cdot (A\tau _e)^{-1}\in M_n(\rat{E})$
and $(A^{-1}\tau _e)=(A\tau _e)^{-1}\in M_n(\rat{E})$ for some
$n\ge 1$, so the result follows. Of course a similar argument
shows that $\rat{E}$ is closed under all the right transductions
$\tilde{\delta}_e$.

Let $X\subseteq E^0$ be the set of vertices which are not sources.
Given a vertex $i\in X$, consider the following element:
$$q_i=p_i-\sum_{e\in r^{-1}(i)}\ol{e}e\in R\left<\ol{E};\tau,\delta\right> .$$

\begin{lem} The elements
 $q_i$ defined above are pairwise orthogonal, nonzero idempotents and
$q_i\leqslant p_i$ for all $i\in X$.
\end{lem}
\begin{proof}
Using the relations $e\ol{f}=\delta _{e,f}p_{s(e)}$ and the
relations in $\cam{E}$ and $\cam{\overline{E}}$ we have that
$$
q_ i^2=p_i^2-\sum_{e\in r^{-1}(i)}\ol{e}ep_i-\sum_{e\in r^{-1}(i)}
p_i\ol{e}e+\left(\sum_{e\in r^{-1}(i)}\ol{e}e\right)^2=q_i,
$$
moreover, $q_ip_i=p_iq_i=q_i$ so that $q_i$ are idempotent
elements and $q_i\leqslant p_i$. Since the $p_i$'s are pairwise
orthogonal, it is clear that the $q_i$'s are also orthogonal. It
follows from  Proposition \ref{prop:anellE}\eqref{item:proj} that
$q_i\ne 0$ for all $i\in X$.
\end{proof}

\begin{nota}
We write $q=\sum_{i\in X}q_i=\sum_{i\in X}p_i-\sum_{e\in
E^1}\ol{e}e\,\,$ which is, by the above lemma, an idempotent.
\end{nota}

\begin{lem}\label{lema:idempip} Let $R$ be a
$\cb$-subalgebra of $\ser{E}$ containing $\cam{E}$, closed under
inversion and closed under all the left transductions
$\delta_{e}$. Set $S=R\langle \ol{E};\tau,\delta \rangle $ and
$I=SqS$, the two-sided ideal generated by the idempotent $q$. Then
the following properties hold:
\begin{enumerate}
\item If $r\in R\setminus \{ 0 \}$ then there exists $y\in S$ such
that $p_iry=p_i$ for some $i\in E^0$. \item If $s \in S\setminus
I$ then there are $s_1,s_2\in S$ such that $s_1ss_2=p_i$ for some
$i\in E^0$. \item If $s\in I\setminus \{ 0 \}$ then there exist
$s_1,s_2\in S$ such that $s_1ss_2=q_i$ for some $i\in X$.
\end{enumerate}
\end{lem}

\begin{proof}
(1) Take $r\in R\setminus \{ 0\}$, with order $k$. Let $w\in E^*$
be a path of length $k$ in the support of $r$ and put $i=s(w)$.
Then $r\ol{w}= \lambda p_i+r'$ where $\lambda \in K\setminus
\{0\}$ and $r'$ is an element in $R$ of order different from $0$.
Thus it follows that $r\ol{w}+(1-p_i)$ is an invertible element in
$R$. Let $t$ be the inverse of $r\ol{w}+(1-p_i)$, and observe that
$$p_ir(\ol{w}t)=p_i,$$
as wanted.

(2) Let $s\in S\setminus I$. By Proposition \ref{prop:anellE} we
know that $s$ can be written as a (finite) right $R$-linear
combination $s=\sum _{\gamma \in E^*}\ol{\gamma} a_{\gamma}$,
where $a_{\gamma}\in p_{s(\gamma )}R$. Observe that $p_j\ol{e}=0$
for all $e\in E^1$ and all $j\in E^0\setminus X$, so that $p_js\in
R$ for all $j\in E^0\setminus X$. Therefore if there exists $j\in
E^0\setminus X$ such that $p_js\ne 0$ then the result follows from
part (1).

So we can assume that $s=p_Xs$, where $p_X=\sum _{i\in X} p_i$. By
an obvious induction, it is enough to show that there is $e\in
E^1$ such that $es\notin I$. If $es\in I$ for all $e\in E^ 1$ then
we have
$$s=qs+(p_X-q)s=qs+ \sum _{e\in E^1}\ol{e}es\in I, $$
a contradiction with our hypothesis. This shows the result.

(3) Since for all $i,j\in X$, all $e\in E^1$ and all $r\in R$, we
have $q_i\ol{e}=0$ and $rq_j=\aug[j]{r}q_j\in Kq_j$, we see that
$q_iSq_j=\delta_{i,j}q_i\cb\cong\cb$. In particular we see that
$I=\sum _{i\in X}Sq_iS$ and indeed, using the above relations and
Proposition \ref{prop:anellE}, we get that every element $s\in I$
can be uniquely written as a finite sum
\begin{equation*}\label{eq:elemI}\tag{2.3}
s=\sum_{i\in X}\sum_{\{\gamma\in E^*\mid s(\gamma)=i
\}}\ol{\gamma} q_i a_{\gamma},
\end{equation*} where $a_{\gamma} \in p_{s(\gamma
)}R$.

Now if $a_{\gamma}=0$ for every $\gamma\in E^*$ of positive
length, then the result follows from (1). Assume that
$a_{\gamma}\ne 0$ for some $\gamma \in E^*$ of positive length. By
induction, it suffices to show there is $e\in E^1$ such that
$es\ne 0$. Let $w\in E^*$ be a path of maximum length in the
support of $s$ (with respect to the above expression) and let $e$
be the final arrow in $w$, in such a way that
$\ol{w}=\ol{e}\ol{w'}$ for some $w'$. Then
$$es=\sum_{i\in X}\sum_{\{\gamma\in E^*\mid s(\gamma)=i\}}\ol{(\gamma \delta _e)}\cdot  q_i
\cdot a_{\gamma}$$ is a nonzero element in $I$.
\end{proof}

\begin{prop}\label{prop:socNeu}
Let $R$ be a $\cb$-subalgebra of $\ser{E}$ containing $\cam{E}$,
closed under inversion and closed under all the left transductions
$\delta_{e}$. The ring $S=R\left<\ol{E};\tau,\delta\right>$ is a
semiprime ring and  $SqS$ is a direct summand of $\soc S$.
Moreover, $SqS$ and $\soc S$ are both von Neumann regular ideals
of $S$.
\end{prop}

\begin{proof} It will be a convenient notation in this proof to
set $q_i=p_i$ for $i\in E^0\setminus X$. Given $s\in S\setminus
\{0\}$, we have by Lemma \ref{lema:idempip} that $q_i\in SsS$ for
some $i\in E^0$. Since $q_i$ are nonzero idempotents, we get that
$(SsS)^2\neq \{0\}$, which shows that $S$ is a semiprime ring.

As observed in the proof of Lemma \ref{lema:idempip}, we have
$q_iSq_j=\delta_{i,j}q_i\cb\cong\delta _{i,j}\cb$ for all $i,j$.
In particular $q_iSq_i$ is a division ring (indeed a field) and so
by \cite[Proposition 21.16(2)]{Lam1} the right ideals $q_iS$ are
minimal. Hence we get $Sq_iS\subseteq\soc{S}$ for all $i\in E^0$.

Now we show that $\soc{S}\subseteq \oplus _{i\in E^0} Sq_iS$. By
definition $\soc{S}$ is the sum of all the minimal right (or left)
ideals of $S$ (see \cite[page 186]{Lam1}). Since $S$ is semiprime,
every minimal right ideal of $S$ is of the form $eS$, where $e$ is
a (nonzero) idempotent in $S$ (see e.g. \cite[Corollary
10.23]{Lam1}). If $e$ is an idempotent such that $eS$ is a minimal
right ideal, then by Lemma \ref{lema:idempip} there exist
$s_1,s_2\in S$ such that $s_1es_2=q_i$ for some $i\in E^0$. Since
$eS$ is a minimal right ideal, we have that $(es_2)S=eS$ and, so
there exists $s_3\in S$ such that $es_2s_3=e$. Moreover by
\cite[Lemma 11.9]{Lam1} we have that $S(es_2)$ is a minimal left
ideal and, since $q_i\in S(es_2)$, we get $S(es_2)=Sq_i$ so that
there exists $s_4\in S$ such that $s_4q_i=es_2$. Finally,
$e=es_2s_3=s_4q_is_3\in Sq_iS$ which proves the desired inclusion.
We also see now that $SqS$ is a direct summand of $\soc{S}$.

Now we show that $\soc{S}$ and $SqS$ are (von Neumann) regular
ideals. Observe that $SqS$ is the orthogonal sum of the ideals
$Sq_iS$, $i=1,\dots, d$, and that these ideals are simple rings
(possibly without unit) and contain a minimal one-sided ideal.
Thus the result follows from Litoff's Theorem (see \cite{FU}).
\end{proof}

We recall the following result from \cite{AGP}, which will be very
useful later on.

\begin{lem}\cite[Lemma 5.3]{AGP}\label{lema:5.3AGP}
Let $A$ be a left semihereditary ring and let $B=\Sigma ^{-1}A$ be
a universal localization of $A$. Suppose that for every finitely
presented right $A$-module $M$ such that $\Hom_A(M,A)=0$, we have
that $M\otimes_A B=0$. Then $B$ is a von Neumann regular ring and
every finitely generated projective $B$-module is induced by a
finitely generated projective $A$-module.
\end{lem}

Let $R$ be a $\cb$-subalgebra of $\ser{E}$ containing $\cam{E}$.
As before, let $X=E^0\setminus \text{Sour}(E)$ be the set of
vertices which are not sources in $E$. For $i\in X$ put
$r^{-1}(i)=\{e^{i}_1,\dotsc,e^{i}_{n_i}\}$ and consider the right
$R$-module homomorphisms
\begin{align*}
\mu_i\colon p_iR&\longrightarrow\bigoplus_{j=1}^{n_i}p_{s(e^i_j)}R\\
r&\longmapsto \left(e^i_1r,\dotsc, e^i_{n_i}r\right).
\end{align*}
Write $\Sigma_1=\{\mu_i\mid i\in X\}$. Observe that the elements
of $\Sigma_1$ are homomorphisms between finitely generated
projective right $R$-modules, so that we can consider the
universal localization $\Sigma_1^{-1}R$.

\begin{prop}\label{prop:Tlocu}
Let $R$ be a $\cb$-subalgebra of $\ser{E}$ containing $\cam{E}$
and closed under the left transductions $\delta_{e}$. Set
$S=R\left<\ol{E};\tau,\delta\right>$, let $I$ be the ideal of $S$
generated by $q$ and let $\Sigma _1$ be as above. Then $\Sigma
_1^{-1}R\cong S/I$.
\end{prop}
\begin{proof}
Set $T:=S/I$. Given $s\in S$ we will denote by $\widetilde{s}$ its
class in $T$. As before, we will identify $R$ with a subring of
$S$. Let $f\colon R\to T$ be the composition of the inclusion
$\iota\colon R\hookrightarrow S$ with the canonical projection
$\pi\colon S\to T$. By the above identification, we can write
$f(r)=\widetilde{r}$.

We want to see that $f$ is a universal $\Sigma_1$-inverting
homomorphism. Define, for $i\in X$, the following right $T$-module
homomorphisms
\begin{align*}
\widetilde{\mu}_i\colon\left(\bigoplus_{j=1}^{n_i}p_{s(e^i_j)}R\right)\otimes_R T&\longrightarrow p_iR\otimes_R T\\
(r_1,\dotsc,r_{n_i})\otimes t &\longmapsto
p_i\otimes\left(\sum_{j=1}^{n_i}\widetilde{\ol{e^{i}_j}r_j}t\right).
\end{align*}
Now, the relations
$\widetilde{p}_i=\sum_{j=1}^{n_i}\widetilde{\ol{e^{i}_j}e^i_j}$
and $\widetilde{e\ol{f}}=\delta_{e,f}\widetilde{p}_{s(e)}$ in $T$
give that
$\widetilde{\mu}_i=\left(\mu_i\otimes\id[T]\right)^{-1}$.
Therefore $f$ is $\Sigma _1$-inverting.

To show that $f$ is universal $\Sigma _1$-inverting, consider a
$\cb$-algebra $A$ and a $\Sigma _1$-inverting algebra homomorphism
$g\colon R\to A$. For $i\in X$ the following diagram is
commutative:
$$
\begin{CD}
\left(\bigoplus_{j=1}^{n_i}p_{s(e^i_j)}R\right)\otimes_R A @>\left(\mu_i\otimes\id[A]\right)^{-1}>> p_iR\otimes_R A\\
   @V\cong VV                                                                                            @VV\cong V\\
\bigoplus_{j=1}^{n_i}g(p_{s(e^i_j)})A @>\left(a^i_1,\dotsc,\,
a^i_{n_i}\right)\cdot>> g(p_i)A
\end{CD}
$$
for some $a^i_j\in g(p_i)Ag(p_{s(e^i_j)})$. From this we conclude
that the compositions
\begin{gather*}
\left(a^i_1,\dotsc,a^i_{n_i}\right)\begin{pmatrix}g(e^i_1)\\\vdots\\g(e^i_{n_i})
\end{pmatrix}=g(p_i)\quad\text{and}\label{equ:inv1}\tag{2.4}\\
\begin{pmatrix}g(e^i_1)\\\vdots\\g(e^i_{n_i})\end{pmatrix}
\left(a^i_1,\dotsc,a^i_{n_i}\right)=\diag\left(g(p_{s(e^i_1)}),
\dotsc,g(p_{s(e^i_{n_i})})\right),\label{equ:inv2}\tag{2.5}
\end{gather*}
give the identities on $g(p_i)A$ and on
$\bigoplus_{j=1}^{n_i}g(p_{s(e^i_j)})A$, respectively.

Take $e\in E^1$; we have that $e=e^i_j$ for some $e^i_j\in
r(i)^{-1}$, where $i=r(e)$. Putting $t_{\ol{e}}=a^i_j$ we conclude
from \eqref{equ:inv2} that $g(e)t_{\ol{e}}=g(p_{s(e)})$ and
$g(e^i_k)t_{\ol{e}}=0$ for $k\neq j$. Moreover, if we take $f\in
E^1$ such that $r(f)\neq i$ we have that
$g(f)t_{\ol{e}}=g(f)g(p_{r(f)})g(p_i)t_{\ol{e}}=0$. We are thus in
the hypothesis of Proposition \ref{prop:panellE} and there exists
a unique algebra homomorphism $\overline{g}\colon S\to A$
extending $g$ and such that $\overline{g}(\ol{e})=t_{\ol{e}}$ for
all $e\in E^1$. From \eqref{equ:inv1}, we get that
$g(p_i)=\sum_{j=1}^{n_i}a^i_jg(e^i_j)$ in $A$, which entails that
$p_i-\sum_{j=1}^{n_i}\ol{e^{i}_j}e^i_j\in\ker(\overline{g})$.
Hence $\overline{g}$ factorizes uniquely through $T$ and we have
$h\colon T\to A$ such that $h\circ\pi=\overline{g}$. Now,
$g=\overline{g}\circ\iota=h\circ\pi\circ\iota=h\circ f$ and $h$ is
unique by uniqueness of inverses and the fact that $T$ is
generated by $R$ and $\ol{E}^1$. We have seen that $f$ is
universal $\Sigma _1$-inverting. Therefore $\Sigma _1^{-1}R\cong
T$.
\end{proof}

\begin{remark}
The uniqueness of the expression in \eqref{eq:uniqueform} for
elements in $S=R\langle \ol{E};\tau ,\delta \rangle$ and the
uniqueness of the expression in \eqref{eq:elemI} for elements in
$I=SqS$ give that the natural maps $R\to T$ and $\cam{\ol{E}}\to
T$ are both injective, where $T=S/I=\Sigma _1^{-1}R$. They also
give that, for two $K$-algebras $R_1$ and $ R_2$ such that
$P(E)\subseteq R_1\subseteq R_2\subseteq P((E))$ and $R_1$ and
$R_2$ are closed under all the transductions $\delta _e$, we have
that $I_1=I_2\cap S_1$, where $S_i=R_i\langle \ol{E}; \tau ,\delta
\rangle  $ and $I_i=S_iqS_i$ for $i=1,2$. It follows that the
natural map $T_1=S_1/I_1\to T_2=S_2/I_2$ is injective.
\end{remark}

\begin{prop}\label{prop:Rsemih}
Let $R$ be a $\cb$-subalgebra of $\ser{E}$ containing $\cam{E}$
and closed under inversion and under all the right transductions
$\tilde{\delta}_{e}$. Then $R$ is left semihereditary.
\end{prop}
\begin{proof}
By \cite{Small} (see also \cite[Proposition 7.63]{Lam2}), it is
enough to show that for every $n\geqslant 1$, $M_n(R)$ is left
Rickart, that is $\lanu[M_n(R)]{A}$ is generated by an idempotent
for all $A\in M_n(R)$. Let $A\in M_n(R)$. We will use the
following notation for the annihilator: $I_A=\lanu[M_n(R)]{A}$.
Given $U\in\gl[n]{R}$ we have that $I_{UA}=\{XU^{-1}\mid X\in
I_{A}\}$.

Observe that $I_A$ is a regular submodule of $M_n(R)$ in the sense
of Section 1, that is if $X\in I_A$ and $o(X)>0$ then
$\tilde{\delta} _e (X)\in I_A$ for every $e\in E ^1$. This will be
used later in the proof.

Since $\aug{I_{A}}$ is a left ideal of $M_n(\cb^d)$, which is
semisimple, there exists an idempotent $D\in\idem{M_n(\cb^d)}$
such that $\aug{I_{A}}=M_n(\cb^d)D$. We can therefore take $B\in
I_{A}$ such that $\aug{B}=D$. Thus we have that $B=D-B'$ for some
$B'\in M_n(R)$ such that $\ord{B'}>0$. Moreover, we can assume
that $DB'=B'$. If we set $U=\id[n]-B'$, since  $R$ is closed under
inversion, we get from Observation \ref{rateqdiv2} that
$U\in\gl[n]{R}$. Now, we see that $BU^{-1}=D$, so that $D\in
I_{UA}$. Since $\aug{U}=\id[n]$ we also have that
$\aug{I_{UA}}=M_n(\cb^d)D$.

We now show that $I_{UA}=M_n(R)D$. We have shown before that
$M_n(R)D\subseteq I_{UA}$. Suppose that there is $X\in
I_{UA}\setminus M_n(R)D$. Substituting $X$ by $X-XD$ we can assume
that $X=X(\id[n]-D)$. Writing $X=\sum_{\alpha\in E^*}\alpha
\lambda_{\alpha}$ for some $\lambda_\alpha\in
p_{r(\alpha)}M_n(\cb^d)$, we get
$\lambda_{\alpha}=\lambda_\alpha(\id[n]-D)$. On the other hand, if
$\ord{X}=m$ we have that $X=\sum_{\alpha\in E^m}\alpha\cdot
\tilde{\delta} _{\alpha}(X)$. Since $I_{UA}$ is a regular
submodule of $M_n(R)$, for every $\alpha\in E^m$ we have that
$\tilde{\delta} _{\alpha}(X)\in I_{UA}$. Take $\alpha\in E^m$ such
that $\lambda_\alpha\neq 0$. We have that $\aug{\tilde{\delta}
_{\alpha}(X)}=\lambda_\alpha$ but, since
$\lambda_{\alpha}(\id[n]-D)=\lambda_\alpha$, this leads us to a
contradiction with $\aug{I_{UA}}=M_n(\cb^d)D$. We have shown that
$I_{UA}=M_n(R)D$ and we deduce that $I_{A}=M_n(R)H$ where
$H=U^{-1}DU\in\idem{M_n(R)}$. It follows that $M_n(R)$ is left
Rickart, as desired.
\end{proof}

\begin{theor}
\label{regular} Let $E$ be a finite quiver and let $R$ be a
$\cb$-subalgebra of $\ser{E}$ containing $\cam{E}$ and closed
under inversion and under all the transductions $\delta_{e}$ and
$\tilde{\delta}_e$. Set $S=R\left<\ol{E};\tau,\delta\right>$,
$I=SqS$ and $T=S/I$. Then $T$ and $S$ are von Neumann regular.
\end{theor}
\begin{proof}
By Proposition \ref{prop:Tlocu} we have that $T$ is a universal
localization of $R$ and, moreover, by Proposition
\ref{prop:Rsemih} $R$ is left semihereditary. Let $M$ be a
finitely presented right $R$-module such that $\Hom_R(M,R)=0$. We
want to show that $M\otimes_R T=0$. Consider the following
presentation of $M$:
\begin{equation*}\label{eq:presentacio}\tag{2.6}
\expesq{R}{s}\xrightarrow{\lpr{A}}\expesq{R}{t}\longrightarrow M\longrightarrow 0,
\end{equation*}
where $A\in\mat[t\times s]{R}$. Adding some zero columns to $A$,
we can assume that $t\leqslant s$. Applying the functor
$\Hom_R(-,R)$ to \eqref{eq:presentacio} we obtain the exact
sequence:
$$
0\longrightarrow\Hom_R(M,R)\longrightarrow R^t\xrightarrow{\rpr{A}}R^s
$$
and, since $\Hom_R(M,R)=0$, we have that $\rpr{A}$ is a
monomorphism. By the right exactness of the functor $-\otimes_R
T$, applied to \eqref{eq:presentacio}, we get the exact sequence:
$$
\expesq{T}{s\,}\xrightarrow{\lpr{A}}\expesq{T}{t\,}\longrightarrow M\otimes_R T\longrightarrow 0.
$$
We want to see that $A\expesq{T}{s\,}=\expesq{T}{t\,}$, that is,
that the columns of $A$ generate $\expesq{T}{t\,}$ as a right
$T$-module.

By a standard argument of linear algebra we know that there are
matrices $P\in\gl[t]{\cb^d}$ and $Q\in\gl[s]{\cb^d}$ such that
$P\aug{A}Q=D$, where
\begin{equation*}\label{eq:PAQ}\tag{2.7}
D=\sum_{i=1}^dp_i\begin{pmatrix}\id[r_i]&0\\0&0\end{pmatrix}\quad\text{with
$r_1,\dotsc,r_d\leqslant t$}.
\end{equation*}
Hence, $PAQ=D-X$, where $X\in\mat[t\times s]{R}$ with $\ord{X}>0$.
Since $t\leqslant s$, we have that $D=(D'\,0)$ with
$D'\in\mat[t]{R}$. Observe that, in the case where $D'=\id[t]$, it
follows from Observation \ref{rateqdiv2} that $PAQ$ is right
invertible over $R$ so that $(PAQ)T^s=T^t$ and we are done.

Otherwise, consider the matrix
$\left(\begin{smallmatrix}X\\0\end{smallmatrix}\right)\in\mat[s]{R}$.
By Observation \ref{rateqdiv2} we know that
$Q'=\id[s]-\left(\begin{smallmatrix}X\\0\end{smallmatrix}\right)\in\gl[s]{R}$.
Therefore,
$$
PAQ(Q')^{-1}=(D-X)\left(\id[s]+\left(\begin{smallmatrix}X\\0\end{smallmatrix}\right)+
\left(\begin{smallmatrix}X\\0\end{smallmatrix}\right)^2+\dotsb\right)=D-(X_1\,X_2),
$$
where $X_1\in\mat[t]{R}$ satisfies that $(\id[t]-D')X_1=X_1$.

Again by Observation \ref{rateqdiv2} we have that
$\id[t]-X_1\in\gl[t]{R}$. Now, setting
$$
A'=(\id[t]-X_1)^{-1}PAQ(Q')^{-1}=D-(X_3\,X_4),
$$
we have that $X_3(\id[t]-D')=(\id[t]-D')X_3=X_3$.

We distinguish two cases, depending on whether $X_3$ is zero or
not. If $X_3\neq 0$ we take $\alpha$ of minimal length amongst the
monomials in the support of the entries of $X_3$. Suppose that
$\alpha$ belongs to the $i$-th column of $X_3$. Consider the
column $v=(0,\dotsc,0,\ol{\alpha},0,\dots,0)^t\in\expesq{S}{s}$,
where $\ol{\alpha}$ is in the $i$-th position. From the condition
$X_3(\id[t]-D')=X_3$ we deduce that $Dv=0$ and thus
$A'v\in\expesq{R}{t}$. In addition we have $(\id[t]-D')A'v=A'v$
and $p_{s(\alpha)}\in\supp{A'v}$.

If $X_3=0$ then by multiplying $A'$ on the right by the matrix
$\left(\begin{smallmatrix}\id[t]&-X_4\\0&\id[s-t]\end{smallmatrix}\right)\in\gl[s]{R}$
we can assume that $X_4$ satisfies that $(\id[t]-D')X_4=X_4$.
Since $\rpr{A}$ is injective we have that $X_4\neq 0$.  As before,
we take $\alpha$ of minimal length amongst all the monomials in
the support of the entries of $X_4$ and we get a column
$v=(0,\dotsc,0,\ol{\alpha},0,\dots,0)^t\in\expesq{S}{s}$
satisfying that $A'v\in\expesq{R}{t}$, $(\id[t]-D')A'v=A'v$ and
$p_{s(\alpha)}\in\supp{A'v}$.

In each of the above two cases, we consider the matrix
$A''=(A'\,\, A'v)\in\mat[t\times (s+1)]{R}$. We have
$$\aug{A''}=(D'\,\,\,0\,\,\, \aug{A'v})$$
so from the conditions $(\id[t]-D')A'v=A'v$ and $p_{s(\alpha)}\in
\supp{A'v}$, we infer that $\rang{K^d}{\aug{A''}}
>\rang{K^d}{\aug{A}}$.
We have the following commutative diagram with exact rows:
$$
\begin{CD}
\expesq{R}{s}@>\lpr{A}>>\expesq{R}{t}@>>>M@>>>0\\
@V\cong VV              @VV\cong V         @|\\
\expesq{R}{s}@>\lpr{A'}>>\expesq{R}{t}@>>>M@>>>0\\
@ViVV               @|         @VVfV\\
\expesq{R}{s}\oplus R@>\lpr{A''}>>\expesq{R}{t}@>>>M'@>>>0
\end{CD}
$$
where $i$ denotes the inclusion in the first component and $f$
exists by the universal property of the cokernel. Applying the
functor $-\otimes_R T$ to the above diagram we get the following
commutative diagram with exact rows:
$$
\begin{CD}
\expesq{T}{s\,}@>\lpr{A}>>\expesq{T}{t\,}@>>>M\otimes_RT@>>>0\\
@V\cong VV @VV\cong V @| \\
\expesq{T}{s\,}@>\lpr{A'}>>\expesq{T}{t\,}@>>>M\otimes_RT@>>>0\\
@VVV @| @VVV\\
\expesq{T}{s+1\,}@>\lpr{A''}>>\expesq{T}{t\,}@>>>M'\otimes_RT@>>>0.
\end{CD}
$$
Since $A'\expesq{T}{s\,}=A''\expesq{T}{s+1\,}$ we have that
$M\otimes_RT\cong M'\otimes_RT$. Moreover it is clear that
$\rpr{A''}\colon R^t\to R^{s+1}$ is injective, so that we can
repeat the above argument. After a finite number of steps we will
arrive at a matrix $B\in\mat[t\times (s+\ell)]{R}$ such that
$\aug{B}=(\id[t]\,0)$ and we see that
$M\otimes_RT\cong\coker{\lpr{B}}=0$.

Now it follows from Lemma \ref{lema:5.3AGP} that $T\cong S/I\cong
\Sigma _1^{-1}R$ is a von Neumann regular ring. Since, by
Proposition \ref{prop:socNeu}, $I$ is a von Neumann regular ideal
of $S$, we get from \cite[Lemma 1.3]{vnrr} that $S$ is von Neumann
regular too.
\end{proof}

\section{The structure of finitely generated projective modules}

Let $E$ be a finite quiver and let $R$ be a $\cb$-subalgebra of
$\ser{E}$ containing $\cam{E}$ and closed under inversion and
under all the transductions $\delta_{e}$ and $\tilde{\delta}_e$.
Set $S=R\left<\ol{E};\tau,\delta\right>$, $I=SqS$ and $T=S/I\cong
\Sigma _1^{-1}R$ (see Section 2).

We have a commutative diagram of inclusion maps
\begin{equation}\notag \begin{CD}
K^d @>>> \cam{E} @>>> R @>>> \ser{E}\\
@VVV @VVV @VVV @VVV\\
\cam{\ol{E}} @>>> L(E) @>>> T @>>> U
\end{CD} \notag \end{equation}
where $U=\Sigma _1^{-1}\ser{E}$, $T=\Sigma _1^{-1}R$ and
$L(E)=\Sigma _1^{-1}\cam{E}$, being $\Sigma _1$ the set of
homomorphisms between finitely generated projective modules
defined in the previous section. In this section we will compute
the structure of the monoid $\mon {T}$, indeed we will show that
the maps in the bottom row of the above diagram induce
isomorphisms $\mon{L(E)}\cong \mon{T}\cong \mon{U}$. The algebras
$L(E)$ are the Leavitt path algebras of \cite{AA1}, \cite{AA2},
\cite{AMP}. The monoid $M_E=\mon{L(E)}$ has been computed in
\cite{AMP}.

Let $F_E$ be the free abelian monoid on the set $E^0$. The nonzero
elements of $F_E$ can be written in a unique form up to
permutation as $\sum _{i=1}^n v_i$, where $v_i\in E^0$. For $v\in
E^0$ such that $r^{-1}(v)\ne \emptyset$, write
$${\bf s}(v):=\sum _{\{e\in E^1\mid r(e)=v\}} s(e)\in F_E .$$
Then the monoid $M_E=\mon{L(E)}$ is isomorphic to $F_E/\sim$,
where $\sim$ is the congruence on $F_E$ generated by all pairs
$(v,{\bf s}(v))$ with $r^{-1}(v)\ne \emptyset $.

\begin{theor}\label{monT}
There is a canonical isomorphism $M_E\cong \mon{T}$.
\end{theor}

\begin{proof}
We have a natural monoid homomorphism $\varphi \colon M_E\cong
\mon{L(E)}\to \mon{T}$. We have seen in the proof of Theorem
\ref{regular} and in Proposition \ref{prop:Rsemih} that the ring
$R$ satisfies the hypothesis of Lemma \ref{lema:5.3AGP}, so we get
that all the finitely generated projective $T$-modules are induced
from finitely generated projective $R$-modules. Observe that
$J(R)=\ker{\varepsilon}$, so that $R/J(R)\cong K^d$ and $R$ is a
semiperfect ring. This follows from the well-known
characterization of $J(R)$ as the set of elements $x$ in $R$ such
that $1-xy$ is invertible for all $y\in R$. Since $R$ is
inversion-closed in $\ser{E}$ all the elements of the form $1-x$,
with $x\in \ker{\varepsilon}$ will be invertible in $R$. It
follows that we have an isomorphism
$$\mon{R}\cong \mon{K^d}.$$
We conclude that all the finitely generated projective $T$-modules
are isomorphic to finite direct sums of modules of the form
$p_iT$, where $p_i$ are the basic idempotents in $K^d$. It follows
that the map $\varphi \colon M_E\to \mon{T}$ is surjective.

Now we will show injectivity of $\varphi$. Assume that
$\oplus_{i=1}^n p_{v(i)}T\cong \oplus _{j=1}^m p_{w(j)}T$, where
$v(i),w(j)\in E^0$. We want to show that $\sum _{i=1}^n v(i)\sim
\sum _{j=1}^m w(j)$ in $F_E$. Since $T$ is von Neumann regular,
the refinement property for f.g. projective modules holds
\cite[Theorem 2.8]{vnrr}, so that we can reduce ourselves to the
case where $n=1$. Let $\alpha \colon p_{v}T\to \oplus _{j=1}^{m}
p_{v(j)}T$ be an isomorphism. Write $\alpha =(\alpha_1,\dots
,\alpha_m) $ with $\alpha _j\in p_{v(j)}Tp_v$. Each $\alpha_j$ can
be written as $\alpha _j=\sum _k w_{jk}\gamma_{jk}$, where
$w_{jk}\in \ol{E}^*$ and $\gamma _{jk}\in R$, all $j,k$, and since
$\alpha _j\in p_{v(j)}T$, we can assume that $\ol{s}(w_{jk})=v(j)$
for all $j,k$.

We proceed by induction on the maximum of the lengths of the paths
$w_{jk}$ appearing in these decompositions. If the maximum is $0$
then the map $\alpha$ is induced from a map $p_vR\to
\oplus_{j=1}^m p_{v(j)}R$ and so the result follows from Lemma
\ref{lemma5.5AGP}. Assume that the maximum $N_0$ of the lengths of
the paths $w_{jk}$ is strictly greater than $0$. Take any path
$w_{j_0k_0}$ of length $N_0$. Since $\ol{s}(w_{j_0k_0})=v(j_0)$,
we see that $v(j_0)$ is not a source in $E$ (i.e. is not a sink in
$\ol{E}$), and so we may consider the right $T$-module isomorphism
\begin{equation} \label{isos}\tag{3.1} (e)_{e\in
r^{-1}(v(j_0))}\colon p_{v(j_0)}T\longrightarrow \bigoplus _{e\in
r^{-1}(v(j_0))}{p_{s(e)}}T, \end{equation} which is given by left
multiplication by the row $(e)_{e\in r^{-1}(v(j_0))}$, with
inverse given by left multiplication by the column
$((\ol{e})_{e\in r^{-1}(v(j_0))})^t$. (The maps $(e)_{e\in
r^{-1}(v(j_0))}$ are the maps $\mu_{v(j_0)}\otimes\id[T]$
considered in Proposition \ref{prop:Tlocu}.) Composing the
isomorphism $\alpha\colon p_vT\to \oplus _{j=1}^m p_{v(j)}T$ with
the isomorphism obtained by applying the above canonical
isomorphisms to all the modules $p_{v(j)}T$ such that there is a
path $w_{jk}$ of length $N_0$ in the corresponding representation
of $\alpha_j$, we obtain a new isomorphism $\alpha '\colon p_vT\to
\oplus _{j=1}^{m'}p_{w(j)}T$ such that the maximum length of the
paths in $\ol{E}$ appearing in the representations of the elements
$\alpha _j'\in T$ is less than $N_0$. Each of the isomorphisms
(\ref{isos}) contributes to a basic transformation $\sum _{i=1}^n
v(i)\rightarrow_1 \sum _{i\ne j}v(i)+{\bf s}(v(j))$ (see
\cite{AMP}), and therefore we have that $\sum _{j=1}^m v(j)\sim
\sum _{j=1}^{m'}w(j)$ in $F_E$. By induction we have that $v\sim
\sum _{j=1}^{m'}w(j)$ and thus $v\sim \sum _{j=1}^{m}v(j)$, which
completes the proof.
\end{proof}

The following lemma completes the proof of Theorem \ref{monT}, its
proof follows the lines of the one of \cite[Lemma 5.5]{AGP}.

\begin{lem}
\label{lemma5.5AGP} Let $p=p_v$ for a fixed $v\in E^0$ and $\alpha
\colon pR\to \oplus _{i=1}^s p_{v(i)}R$ such that $\alpha$ becomes
invertible over $T$. Then $v\sim \sum_{i=1}^s v_i$ in $F_E$.
\end{lem}

\begin{proof}
Write $\alpha =(\alpha _1,\dots ,\alpha _s)^t$, where each $\alpha
_i\in p_{v(i)}Rp$. We will construct , by induction on $i$, paths
$w_i\in \ol{E}^*$ and invertible elements $g_i\in
p_{v(i)}Rp_{v(i)}$ such that the following statements hold:

$(A_i)$ There exists an invertible map $\alpha ^{(i)}\colon pT\to
\oplus _{i=1}^s p_{v(i)}T$ satisfying the following properties:

(1) $\alpha ^{(i)}_{i+1},\dots ,\alpha ^{(i)}_s\in R$.

(2) The inverse of $\alpha ^{(i)}$ is the row $(w_1g_1,\dots
,w_ig_i,\beta _{i+1},\dots , \beta _s)$ for some elements $\beta
_{\ell}\in pTp_{v(\ell)}$, $\ell =i+1,\dots ,s$.

The statement is obvious for $i=0$. Assume that $0\le i<s $ and
that $(A_i)$ holds. We will prove $(A_{i+1})$. Without loss of
generality, we can assume that the order of the series $\alpha
^{(i)}_{i+1}$ is less than or equal to the order of $\alpha
^{(i)}_{i+t}$ for all $t\ge 2$. Choose a path $w_{i+1}\in
\ol{E}^*$ with length equal to the order of $\alpha_{i+1}^{(i)}$
such that $w_{i+1}=pw_{i+1}p_{v(i+1)}$ and such that
$\alpha_{i+1}^{(i)}w_{i+1}$ is invertible in
$p_{v(i+1)}Rp_{v(i+1)}$. Let $g_{i+1}\in p_{v(i+1)}Rp_{v(i+1)}$ be
the inverse of $\alpha _{i+1}^{(i)}w_{i+1}$ and note that
$$p_{v(i+1)}=\alpha _{i+1}^{(i)}w_{i+1}g_{i+1}=\alpha _{i+1}^{(i)}\beta _{i+1}.$$
It follows that $$u:=\beta _{i+1}\alpha
_{i+1}^{(i)}+(p-w_{i+1}g_{i+1}\alpha _{i+1}^{(i)})$$ is invertible
in $pTp$ with inverse $$u^{-1}= w_{i+1}g_{i+1}\alpha
_{i+1}^{(i)}+(p-\beta _{i+1}\alpha _{i+1}^{(i)}) .$$ Therefore,
$\alpha ^{(i+1)}:=\alpha ^{(i)}u$ is invertible with inverse
$$u^{-1}(w_1g_1,\dots ,w_ig_i,\beta _{i+1},\dots ,\beta _s).$$
Note that, for $t>1$, we have
$$\alpha_{i+t}^{(i)}u=\alpha _{i+t}^{(i)}(p-w_{i+1}g_{i+1}\alpha _{i+1}^{(i)}) .$$
Since the order of $\alpha _{i+t}^{(i)}$ is greater than or equal
to the length of $w_{i+1}$, we conclude that $\alpha
_{i+t}^{(i+1)}\in R$, and condition (1) of $(A_{i+1})$ holds. On
the other hand, for $m\le i$ we have $\alpha _{i+1}^{(i)}w_m=0$
and so $u^{-1}w_mg_m=w_mg_m$. We also have
$$u^{-1}\beta _{i+1}=w_{i+1}g_{i+1}\alpha _ {i+1}^{(i)}\beta _{i+1}
+(p-\beta _{i+1}\alpha _{i+1}^{(i)})\beta _{i+1}=w_{i+1}g_{i+1}
,$$ and so condition (2) of $(A_{i+1})$ is also satisfied.
Therefore, the induction works.

Take $h_i=g_i\alpha _i^{(s)}\in p_{v(i)}Tp$ for $i=1,\dots ,s$.
Then
\begin{equation}
\label{(1)}\tag{3.2} \sum _{i=1}^s w_ih_i=p ,
\end{equation}
$h_iw_i\ne 0$ for all $i$, and $h_iw_j=0$ for $i\ne j$. We claim
that these conditions imply $v\sim \sum _{i=1}^s v(i)$ in $F_E$.
We proceed by induction on the maximum of the lengths of the
$w_i$. If this maximum is $0$ then $s=1$ and $h_1=p$. So assume
that either $s>1$ or $s=1$ and the length of $w_1$ is $\ge 1$. In
either case, all $w_i$ are different from $p$. Note that
$w_i=\ol{\gamma_i}$ for a path $\gamma _i$ in $E$ of length $\ge
1$ such that $s(\gamma _i)=v(i)$ and $r(\gamma _i)=v$. Let
$e(i)\in E^1$ be the ending arrow of the path $\gamma _i$, so that
$r(e(i))=v$. For each $e\in E^1$ such that $r(e)=v$, define
$$A_e:=\{ i\in \{1,\dots ,s\} \mid e(i)=e \}=\{ i\in \{1,\dots ,s\} \mid ew_i\ne 0 \}.$$
Then the set $\{ 1,\dots , s\}$ is the disjoint union of the sets
$A_e$, for $e\in r^{-1}(v)$.

Fix an arrow $e\in E^1$ such that $r(e)=v$. Left multiplying
\eqref{(1)} by $e$ and right multiplying it by $\ol{e}$, we get
$$\sum _{i\in A_e}(ew_i)(h_i\ol{e})=ep\ol{e}=p_{s(e)}.$$
Observe also that for $i,j\in A_e$ we have
$(h_i\ol{e})(ew_j)=h_iw_j$. So this term is $0$ if $i\ne j$ and
nonzero if $i=j$. By induction, $p_{s(e)}\sim \sum _{i\in
A_e}p_{v(i)}$. Therefore $$p=p_v\sim \sum _{e\in r^{-1}(v)}
p_{s(e)}\sim \sum _{e\in r^{-1}(v)}\sum _{i\in A_e} p_{v(i)}=\sum
_{i=1}^s p_{v(i)} .$$ This shows the result.
\end{proof}

\section{The regular algebra of a quiver.}

In this section we will observe that our construction can be made
functorial in $E$ if we choose a suitable class of morphisms
between quivers. This also enables us to extend the construction
to all the column-finite quivers.

\begin{defi}
For a finite quiver $E$ and a field $K$, we define {\em the
regular algebra} of $E$ as the algebra $Q(E)$ obtained by using
the construction in Section 2 taking as coefficients $R=\rat{E}$.
So we have $$Q(E)=S/SqS=(\Sigma _1)^ {-1}(\rat{E}),$$ where
$S=(\rat{E})\langle \ol{E};\tau, \delta\rangle$.
\end{defi}

The algebra $Q(E)$ fits into a commutative diagram of injective
algebra morphisms:

\begin{equation}\notag \begin{CD}
K^d @>>> \cam{E} @>{\iota _{\Sigma}}>> \rat{E} @>>> \ser{E}\\
@VVV @V{\iota _{\Sigma_1}}VV @V{\iota _{\Sigma _1}}VV @V{\iota _{\Sigma _1}}VV\\
\cam{\ol{E}} @>>> L(E) @>{\iota _{\Sigma}}>> Q(E) @>>> U(E)
\end{CD} \notag \end{equation}
where $U(E)=\Sigma_1^{-1}\ser{E}$, $Q(E)=\Sigma_1^{-1}\rat{E}$ and
$L(E)=\Sigma_1^{-1}\cam{E}$, being $\Sigma _1$ the set of
homomorphisms between finitely generated projective modules
defined in Section 2. We summarize the properties of the algebra
$Q(E)$ in the following theorem:

\begin{theor}
Let $E$ be a finite quiver. Then the regular algebra $Q(E)$ of $E$
is a von Neumann regular hereditary ring, and $Q(E)=(\Sigma\cup
\Sigma _1)^{-1}P(E)$ is a universal localization of the path
algebra $P(E)$. Moreover we have $\mon{Q(E)}\cong M_E$
canonically.
\end{theor}

\begin{proof}
By Theorem \ref{regular} we have that $Q(E)$ is von Neumann
regular, and by Theorem \ref{monT} we have that $\mon{Q(E)}\cong
M_E $ canonically. Using Theorem \ref{ratuniversal}, we get
$Q(E)=(\Sigma\cup \Sigma _1)^{-1}P(E)$, and, since $P(E)$ is a
hereditary ring, a result of Bergman and Dicks \cite{BD} gives
that $Q(E)$ is hereditary too.
\end{proof}

A quiver $E$ is said to be {\em column-finite} in case each vertex
in $E$ receives only a finite number of arrows, that is
$r_E^{-1}(v)$ is finite for all $v\in E^0$. For a column-finite
quiver $E$, one can define the Leavitt path algebra $L(E)$
(\cite{AA1}, \cite{AMP}) and also the quiver monoid $M_E$ just as
in Section 3 (see \cite{AMP}). Note that the Leavitt path algebra
$L(E)$ is unital if and only if the quiver $E$ is finite. As shown
in \cite[Section 2]{AMP} these constructions are functorial with
respect to complete graph homomorphisms, defined below.

Let $f=(f^0,f^ 1)\colon E\to F$ be a graph homomorphism. Then $f$
is said to be {\em complete} if $f^0$ and $f^1$ are injective and
$f^1$ restricts to a bijection between $r_E^{-1}(v)$ and
$r_F^{-1}(f^0(v))$ for every vertex $v\in E^0$ that receives
arrows.

If $f\colon E\to F$ is a complete graph homomorphism between
finite quivers $E$ and $F$, then $f$ induces a non-unital algebra
homomorphism $P(f)\colon P(E)\to P(F)$ between the corresponding
path algebras and a non-unital homomorphism $L(f)\colon L(E)\to
L(F)$ between the corresponding Leavitt path algebras. Note that
the image of the identity under these homomorphisms is the
idempotent
$$p_E:=\sum _{v\in E^0} p_{f^0(v)}\in P(F).$$ We get a morphism
$P(E)\to P(F)\to L(F)\to Q(F)$ such that every map in $\Sigma
_1(E)$ becomes invertible over $Q(F)$.

Observe that we have a commutative diagram

\begin{equation}\notag \begin{CD}
P(E) @>>> P(F) \\
@V{\varepsilon _E}VV @V{\varepsilon _F}VV \\
K^{d_E} @>>> K^{d_F}
\end{CD} \notag \end{equation}
so a matrix $A\in M_n(P(E))$ such that $\varepsilon _E(A)$ is
invertible is sent to a matrix $P(f)(A)\in M_n(p_EP(F)p_E)$ such
that $\varepsilon_F(P(f)(A))$ is invertible over $p_EK^{d_F}p_E$.
It follows that the unital algebra homomorphism $P(E)\to
p_EQ(F)p_E$ factorizes uniquely through $Q(E)=(\Sigma \cup \Sigma
_1)^{-1}P(E)$, so that we get a unital algebra homomorphism
$Q(E)\to p_EQ(F)p_E$, thus a non-unital algebra homomorphism
$Q(f)\colon Q(E)\to Q(F)$ such that $Q(f)(1)=p_E$.

This gives the functoriality property of the regular algebra of a
quiver, for finite quivers. By \cite[Lemma 2.1]{AMP} every
column-finite quiver $E$ is the direct limit, in the category of
quivers with complete graph homomorphisms, of the directed family
$\{E_{\lambda}\}$ of its complete finite subquivers. Thus we get a
directed system $\{ Q(E_{\lambda})\}$ of von Neumann regular
algebras and (non-unital) algebra morphisms, and we define the
regular algebra of $E$ as:
$$Q(E)=\varinjlim Q(E_{\lambda}).$$ Since the functor $\mathcal V$
commutes with direct limits and $M_E\cong \varinjlim
M_{E_{\lambda}}$  (\cite[Lemma 2.4]{AMP}) we get:

\begin{theor}
Let $E$ be any column-finite quiver. Then there is a (possibly
non-unital) von Neumann regular algebra $Q(E)$ such that
$$\mon{Q(E)}\cong M_E.$$
\end{theor}
This solves the realization problem for the monoids associated to
column-finite quivers. Clearly the functoriality of $Q$ extends to
the category of column-finite quivers and complete graph
homomorphisms.


\begin{thebibliography}{99}

\bibitem{AA1} G. {\sc Abrams}, G. {\sc Aranda Pino}, The Leavitt path
algebra of a graph, {\em J. Algebra}, {\bf 293}  (2005), 319--334.

\bibitem{AA2} G. {\sc Abrams}, G. {\sc Aranda Pino}, Purely infinite
simple Leavitt path algebras, {\em J. Pure Appl. Algebra}, to appear.

\bibitem{AB} P. {\sc Ara}, M. {\sc Brustenga}, $K_1$ of corner
skew Laurent polynomial rings and applications, {\em Commun. in
Algebra}, {\bf 33} (2005), 2231--2252.

\bibitem{AGP} P. {\sc Ara}, K. R. {\sc Goodearl}, E. {\sc Pardo},
$K_0$ of purely infinite simple regular rings, {\em K-Theory},
{\bf 26} (2002), 69--100.

\bibitem{agop} P. {\sc Ara}, K. R. {\sc Goodearl}, K. C. {\sc
O'Meara}, E. {\sc Pardo}, Separative cancellation for projective
modules over exchange rings, {\em Israel J. Math.}, {\bf 105}
(1998), 105--137.


\bibitem{AMP} P. {\sc Ara}, M. A. {\sc Moreno}, E. {\sc Pardo},
Nonstable $K$-theory for graph algebras, {\em Algebras Repr.
Theory}, to appear, arXiv: math.RA/0412243.

\bibitem{APS} G. {\sc Aranda Pino}, E. {\sc Pardo}, M. {\sc Siles
Molina}, Exchange Leavitt path algebras and stable rank, {\em J.
Algebra}, to appear.


\bibitem{AuslanderReitenSmalo} {\sc M. Auslander, I. Reiten, S.O. Smal\o.}
``Representation theory of Artin algebras''
Cambridge Studies in Advanced Mathematics, 36. Cambridge
University Press, Cambridge,  1995.


\bibitem{Bergman} {\sc G. M. Bergman}, Coproducts and some
universal ring constructions, {\em Trans. Amer. Math. Soc.}, {\bf
200} (1974), 33--88.

\bibitem{BD} G. M. {\sc Bergman}, W. {\sc Dicks}, Universal
derivations and universal ring constructions, {\em Pacific J.
Math.}, {\bf 79} (1978), 293--337.


\bibitem{BR} {\sc J. Berstel, C. Reutenauer}, ``Rational Series
and Their Languages", EATCS Monographs on Theoretical Computer
Science, vol. 12, Springer-Verlag, Berlin and Heidelberg, 1988.


\bibitem{free} P. M. {\sc Cohn}, ``Free Rings and Their
Relations'', Second Edition, LMS Monographs 19, Academic Press,
London, 1985.


\bibitem{Cohnenc} P. M. {\sc Cohn}, ``Skew fields. Theory of general division rings",
Encyclopedia of Mathematics and its Applications, 57. Cambridge
University Press, Cambridge, 1995.



\bibitem{CD} P. M. {\sc Cohn}, W. {\sc Dicks}, Localization in
semifirs II, {\em J. London Math. Soc.}, {\bf 13} (1976),
411--418.


\bibitem{FU} C. {\sc Faith}, Y. {\sc Utumi}, On a new proof of
Litoff's Theorem, {\em Acta Hung. Math.}, {\bf 14} (1964),
369--371.

\bibitem{vnrr} K. R. {\sc Goodearl}, ``Von Neumann Regular Rings'',
Pitman, London 1979; Second Ed., Krieger, Malabar, Fl., 1991.

\bibitem{directsum}
K. R. {\sc Goodearl}, ``von Neumann regular rings and direct sum
decomposition problems", Abelian groups and modules (Padova,
1994), Math. Appl. 343, 249--255, Kluwer Acad. Publ., Dordrecht,
1995.


\bibitem{Lam1} T.Y. {\sc Lam}, ``A first course in noncommutative
rings", Second edition, Graduate Texts in Mathematics 131,
Springer-Verlag, New York, 2001.

\bibitem{Lam2} T.Y. {\sc Lam}, ``Lectures on modules and rings'', Graduate
Texts in Mathematics 189, Springer-Verlag, New York, 1999.


\bibitem{Luck} W. {\sc L\"{u}ck}, ``$L^2$-invariants: Theory and
Applications to Geometry and $K$-Theory",
Springer Verlag, Berlin,
2002.

\bibitem{Raeburn} I. {\sc Raeburn}, Graph algebras, CBMS Reg. Conf. Ser. Math.,
vol. 103, Amer. Math. Soc., Providence, RI, 2005.

\bibitem{scho} A. H. {\sc Schofield}, ``Representations of Rings
over Skew Fields'', LMS Lecture Notes Series 92, Cambridge Univ.
Press, Cambridge, UK, 1985.

\bibitem{Small} L. {\sc Small}, Semihereditary rings, {\em Bull. Amer. Math.
Soc.}, {\bf 73} (1967), 656--658.

\bibitem{Wehisrael} F. {\sc Wehrung}, Non-measurability properties of interpolation
vector spaces, {\em Israel J. Math.}, {\bf 103} (1998), 177--206.

\bibitem{Wehsemforum}
F. {\sc Wehrung}, Embedding simple commutative monoids into simple
refinement monoids,  {\em Semigroup Forum}, {\bf 56} (1998),
104--129.




\end{thebibliography}
\end{document}